\documentclass[a4paper,10pt]{article}

\usepackage[applemac]{inputenc}   
\usepackage[T1]{fontenc}      
\usepackage{geometry}         
\usepackage[english]{babel}
\usepackage{amsthm}
\usepackage{amsmath}
\usepackage{amsfonts}
\usepackage{amssymb}
\usepackage[all]{xy}
\usepackage[pdftex]{hyperref}
\newtheorem{theorem}{Theorem}[section]
\newtheorem{proposition}[theorem]{Proposition}
\newtheorem{lemma}[theorem]{Lemma}

\newtheorem{corollary}[theorem]{Corollary}
\theoremstyle{definition} \newtheorem{remark}[theorem]{Remark}
\theoremstyle{definition} \newtheorem{example}[theorem]{Example}

\newcommand{\somme}[2]{\displaystyle\sum_{#1}^{#2}}

\def\N {{\mathbb{N}}}

\def\K {{\mathbb{K}}}

\def\S {{\mathbb{S}}}

\def\gg {{\mathfrak{g}}}

\def\conj { {\rhd} }

\def\Sh {{\mathrm{Sh}}}

\title{On the conjectural Leibniz cohomology for groups}
\author{Simon Covez}
\date{}
\begin{document}
\maketitle
\begin{abstract} The goal of this paper is to present results which are consistent with conjectures about \textit{Leibniz (co)homology for discrete groups}, due to J. L. Loday in 2003. We show that \textit{rack cohomology} has properties very close to the properties expected for the conjectural Leibniz cohomology. In particular, we prove the existence of a graded \textit{dendriform} algebra structure on rack cohomology, and we construct a graded associative algebra morphism $H^{\bullet}(-) \to HR^{\bullet}(-)$ from group cohomology to rack cohomology which is injective for $\bullet = 1$.
\end{abstract}
\section*{Introduction}
\paragraph{Problems and results.} The starting point of this article are conjectures by J.-L. Loday in \cite{Lodayconjectural}. There J.-L. Loday conjectures the existence of a Leibniz homology $HL_{\bullet}(-)$ for groups and some properties it should satisfy. Here we focus on two properties
\begin{enumerate}
\item $HL_{\bullet}(G)$ is a graded Zinbiel coalgebra (and, \textit{a fortiori}, a graded cocommutative coalgebra),
\item There exists a morphism of graded cocommutative coalgebras from $HL_{\bullet}(G)$ to the usual group homology $H_{\bullet}(G)$.
\end{enumerate}
Here we work in the cohomological context and we show two results which are consistent with these conjectures \big(Theorem \ref{Theorem : CR(X,A) is a differential graded dendriform algebra} and Theorem \ref{Theorem : a graded associative algebra morphism from group cohomology to rack cohomology}\big). We prove that there exists a cohomology theory $HR^{\bullet}$, defined for groups, and satisfying :
\begin{enumerate}
\item $HR^{\bullet}(G,A)$ is a graded dendriform algebra (and, \textit{a fortiori}, a graded associative algebra),
\item There exists a non zero morphism of graded associative algebras from $H^{\bullet}(G,A)$ to $HR^{\bullet}(G,A)$ which is injective for $\bullet =1$. 
\end{enumerate} 
where $A$ is an associative algebra considered as a trivial $G$-module. To understand what is $HR^{\bullet}(-)$ and where it comes from, let us explain some background information about these conjectures.
\paragraph{Algebraic K-theory and additive K-theory.} The starting point is a program of J.-L. Loday exposed in \cite{Lodayconjectural} which intends to solve two problems arising in algebraic $K$-theory using the knowledge of the solutions of these problems in the additive $K$-theory context. These two problems are :
\begin{enumerate}
\item Have a small presentation of the algebraic $K$-groups of a field, or, at least, a small presentation of a chain complex whose homology would give these $K$-groups. (By small, J.-L. Loday means no matrices in the presentation),
\item Determine the obstructions to the periodicity.
\end{enumerate}
Rationally, the algebraic $K$-theory $K_{\bullet}(A)_{\mathbb{Q}}$ of a unitary ring $A$ can be defined using the group homology of $GL(A)$, as the primitive part of the graded connected cocommutative Hopf algebra $H_{\bullet}(GL(A),\mathbb{Q})$, hence, using the Milnor-Moore theorem, we have the following isomorphism of graded Hopf algebras
$$
H_{\bullet}(GL(A),\mathbb{Q}) \simeq S(K_{\bullet}(A)_{\mathbb{Q}}).
$$
Using this relation, the additive $K$-theory can be seen as the "tangent" version of algebraic $K$-theory. Let $\mathbb{K}$ be a field of characteristic zero and $A$ be a unitary $\mathbb{K}$-algebra. The additive $K$-theory $K^+_{\bullet}(A)$ is defined similarly to algebraic $K$-theory, replacing the group $GL(A)$ by the Lie algebra $\mathfrak{gl}(A)$, and group homology by Lie algebra homology. Therefore $K^+_{\bullet}(A)$ is defined as the primitive part of the graded connected commutative and cocommutative Hopf algebra $H_{\bullet}(\mathfrak{gl}(A),\mathbb{K})$, and by the Milnor-Moore theorem we have the following graded isomorphism
$$
H_{\bullet}(\mathfrak{gl}(A),\mathbb{K}) \simeq S\big(K_{\bullet}^+(A)\big).
$$ 
In the additive $K$-theory context, the solutions of our problems are given by :
\begin{enumerate}
\item The \textit{cyclic homology} $HC_{\bullet}$ by the \textit{Loday-Quillen-Tsygan Theorem} (cf.\cite{LodayQuillen,FeiginTsygan})
\begin{theorem}[Loday-Quillen-Tsygan] Let $\mathbb{K}$ be a field containing $\mathbb{Q}$ and $A$ be a unital associative $\mathbb{K}$-algebra. Then there is a natural isomorphism :
$$
K^+_{n}(A) \simeq HC_{n-1}(A), \, n \geq 1.
$$
\end{theorem}
\item The \textit{Hochschild homology} $HH_{\bullet}$ by the \textit{Connes' periodicity exact sequence} 
$$
\cdots \to HH_n(A) \to HC_n(A) \to HC_{n-2}(A) \to HH_{n-1}(A) \to \cdots .
$$
\end{enumerate}
The following theorem is the starting point of this program. In the same way as the Loday-Quillen-Tsygan Theorem links $HC_{\bullet}(A)$ and $\mathfrak{gl}(A)$, the \textit{Loday-Cuvier Theorem} (cf. \cite{LodayCyclic,CuvierLoday_thm}) gives a link between $HH_{\bullet}(A)$ and $\mathfrak{gl}(A)$.
\begin{theorem}[Loday-Cuvier] For any associative unital algebra $A$ over a characteristic zero field $\mathbb{K}$ there is an isomorphism of graded modules :
$$
HL_{\bullet}(\mathfrak{gl}(A),\mathbb{K}) \simeq T\big(HH_{\bullet-1}(A)\big).
$$
\end{theorem}
\noindent
In this theorem $HL_{\bullet}(-)$ is the Leibniz homology, the homology theory naturally associated to \textit{Leibniz algebras}  (cf \cite{LodayEns,LodayOverview,LodayCyclic}). Leibniz algebras are a non-commutative version of Lie algebras and their homology gives new invariants for Lie algebras. The existence of $HL^{\bullet}(-)$ and its following properties (cf. \cite{Lodaycup}) can explain the reason of the conjectures :
\begin{enumerate}
\item $HL^{\bullet}(\mathfrak{g})$ is a graded Zinbiel algebra,
\item There exists a non zero morphism of graded commutative algebra from $H^{\bullet}(\mathfrak{g})$ to $HL^{\bullet}(\mathfrak{g})$ which is an isomorphism for $\bullet = 1$.
\end{enumerate}
Lie algebra homology being the "tangent" version of group homology, Leibniz algebra homology should be the "tangent" version of the expected homology $HR_{\bullet}(-)$. This hypothesis led J.-L. Loday to the formulation of the \textit{coquecigrue problem} (cf. \cite{LodayEns, Lodayconjectural, Loday_Problems_in_operad}).
\paragraph{The coquecigrue problem.} There are (at least) two ways to construct a Lie algebra from a group. One is to consider the graded abelian group $\oplus_n^{} G^{(n)} / G^{(n+1)}$ associated to the descending central series $\{G^{(n)} = [G,G^{(n-1)}]\}_{n \in \mathbb{N}}$. This object is provided with a Lie algebra structure where the bracket is induced by the commutator in $G$, the Jacobi identity being a consequence of the so-called \textit{Phillip Hall relation}.
\par
Another way is to consider the tangent space at the neutral element of a Lie group. The bracket is induced by the conjugation in the group, and the Leibniz identity by the self-distributivity of the conjugation, that is, the relation :
$$
x \rhd (y \rhd z) = (x \rhd y) \rhd (x \rhd z) \, \text{ where } \, x \rhd y = xyx^{-1}.
$$ 
\par
The coquecigrue problem is to find a generalization of one of these constructions for Leibniz algebra. Therefore the problem is to know if there exists a mathematical object, dubbed coquecigrue, generalizing groups, such that a Leibniz algebras is naturally associated to it (by one of these constructions).
\par
In \cite{CovezIntegration} a local solution to the coquecigrue problem has been given for the second construction using \textit{racks}, especially the links between (Lie) rack cohomology $HR^{\bullet}(-)$ and Leibniz algebra cohomology $HL^{\bullet}(-)$. Originally, racks were defined by knot theorists to construct invariants for knots (cf. \cite{CarterSaito} and references therein), and M. K. Kinyon in \cite{Kinyon} is the first to state a relation between Lie racks and Leibniz algebras.
\par
Since a group is a rack it is natural to suppose that rack cohomology is our conjectural Leibniz cohomology $HR^{\bullet}(-)$, and so, to study the existence of an algebraic product structure on it and its link with usual group cohomology.
\vskip 0.2 cm
The plan for this article is the following :
\vskip 0.2 cm
\paragraph{Section 1 : Dendriform and Zinbiel algebras.} This section is based on \cite{LodayDialgebras}. We recall the definitions of dendriform and Zinbiel algebras and give some examples. We recall also the relations between these algebras and others types of algebras such as \textit{associative, commutative, Lie, Leibniz and diassociative} algebras. 
\paragraph{Section 2 : Racks.} Racks are sets provided with a product which encapsulates some properties of the conjugation in a group (especially the \textit{self-distributivity}). This section is based on \cite{FennRourke} for the general theory of racks and on \cite{EtingofGrana} for the definition of the rack cohomology $HR^{\bullet}(-)$. We recall some basic definitions about racks and their cohomology. 
\paragraph{Section 3 : The category of trunks.} \textit{Trunks} have been introduced by R. Fenn, C. Rourke and B. Sanderson in \cite{FennRourkeSanderson} in order to define a classifying space for racks. Trunks are objects loosely analoguous to categories. Categories can be seen as graphs with a fixed set of triangle (the graph of the composition) satisfying the associativity condition. 
$$
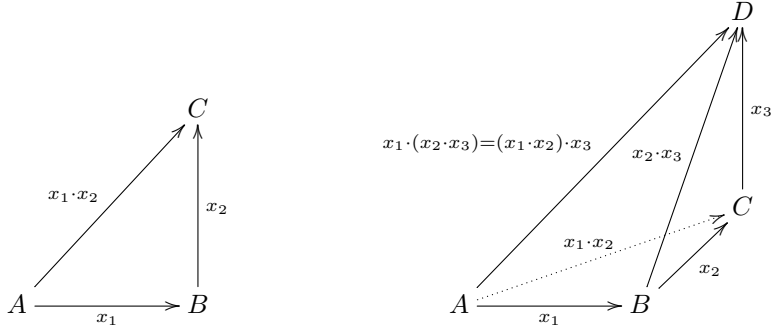
\xymatrix{
& & & & & & & & D
\\
& & C & & & & & &
\\
& & & & & & & & C \ar[uu]_{x_3}
\\
A \ar[rr]_{x_1}  \ar[rruu]^{x_1 \cdot x_2} & & B  \ar[uu]_{x_2} & & & A \ar[rr]_{x_1}  \ar@{.>}[rrru]^{x_1\cdot x_2}  \ar[rrruuu]^{x_1\cdot(x_2\cdot x_3)=(x_1\cdot x_2)\cdot x_3} & & B \ar[ruuu]^{x_2 \cdot x_3}  \ar[ru]_{x_2} &
}
$$
Like categories, trunks are graphs with a fixed set of \textit{cubes}
$$
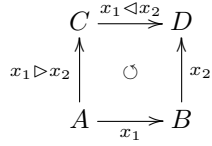
\xymatrix{
C \ar[r]^{x_1 \lhd x_2} \ar@{}[rd]|{\circlearrowleft} & D 
\\
A \ar[r]_{x_1}  \ar[u]^{x_1 \rhd x_2}  & B  \ar[u]_{x_2}
}
$$
and \textit{corner trunks} are trunks where the set of cubes satisfies the \textit{bidistributivity} condition. 
$$
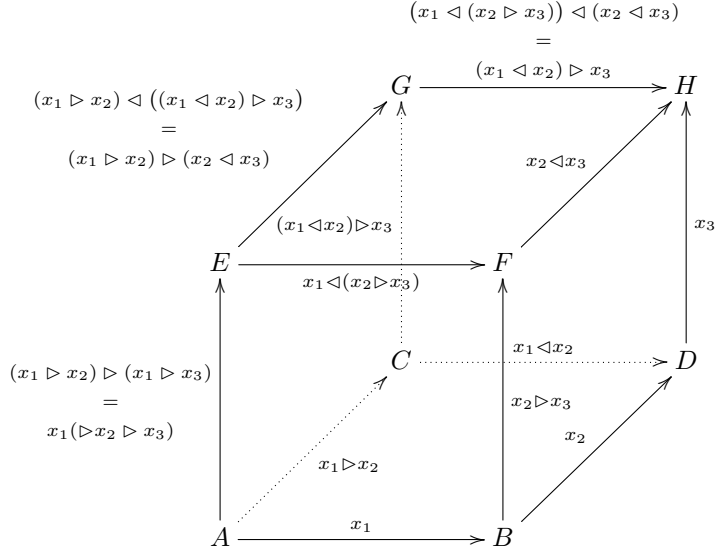
\xymatrix{
& & G \ar[rrr]^{\txt{\scriptsize $\big(x_1 \lhd (x_2 \rhd x_3)\big) \lhd (x_2 \lhd x_3)$ \\ \scriptsize$=$ \\ \scriptsize $ (x_1 \lhd x_2) \rhd x_3$}} & & & H
\\
& & & & &
\\
E \ar[rrr]_{x_1 \lhd (x_2 \rhd x_3)} \ar[rruu]^{\txt{\scriptsize $(x_1 \rhd x_2) \lhd \big((x_1 \lhd x_2) \rhd x_3\big)$ \\ \scriptsize $=$ \\ \scriptsize$(x_1 \rhd x_2) \rhd (x_2 \lhd x_3)$}} &  & & F \ar[urur]^{x_2 \lhd x_3} & &
\\
& & C \ar@{.>}[rrr]^{x_1 \lhd x_2} \ar@{.>}[uuu]^{(x_1 \lhd x_2) \rhd x_3} & & & D \ar[uuu]_{x_3}
\\
& & & & & 
\\
A \ar[rrr]^{x_1}  \ar@{.>}[rruu]_{x_1\rhd x_2}  \ar[uuu]^{\txt{\scriptsize$(x_1 \rhd x_2) \rhd (x_1 \rhd x_3)$ \\ \scriptsize $=$ \\ \scriptsize $x_1(\rhd x_2 \rhd x_3)$}} & & & B \ar[uuu]_{x_2 \rhd x_3}  \ar[rruu]^{x_2} &
}
$$
In the same way as a group can be seen as a category with one object and all its arrows invertible, to any rack $X$ we can associate a corner trunk $X_{\mathbf{Tr}}$ \big(Example \ref{Example : racks and natural numbers are trunks}\big). Moreover, one may construct a (pre)simplicial set from a category (the nerve), in the same way there is a \textit{precubical set} $N(\mathcal{T})$ constructed from a trunk $\mathcal{T}$, that is, a contravariant functor from the category $\square$ to the category $\mathbf{Set}$, by $N(\mathcal{T})(-) := \mathrm{Hom}_{\mathbf{Trunk}}(-,\mathcal{T})$ \big(cf. \eqref{precubical nerve of a trunk}\big). In the case where the trunk is the trunk canonically associated to a rack, the description of the precubical nerve is easy to compute \big(Theorem \ref{Theorem : nerve of a rack}\big), and these result allows us to give a cubical description of the rack cohomology \big(Corollary \ref{Corollary : isomorphism between CR(X,A) and C(N(X),A)}\big). Then we use this description to define a dendriform structure \big(formulas \eqref{definition vdash} and \eqref{definition dashv}\big). Our interest for this cubical definition is that the proof that the differential graded dendriform structure on $HR^{\bullet}(-)$ is well defined becomes a succession of combinatorial lemmas.
\paragraph{Section 4 : A graded dendriform algebra structure on $HR^{\bullet}(X,A)$.} This section and the following are the heart of this paper. We use the cubical description of rack cohomology using trunks described in Section 3 to define a graded dendriform algebra structure on the cochain complex $\{CR^n(X,A),d_R^n\}_{n \in \mathbb{N}}$ computing rack cohomology of a rack $X$ with coefficients in an associative algebra (considered as a trivial $X$-module) (Theorem \ref{Theorem : graded dendriform algebra structure}). Then we show that this structure is compatible with the differential (Theorem \ref{Theroem : compatibility with the differential}).
\par
Dendriform and Zinbiel algebras are closely related to shuffles as we can see examples \ref{Example : Dendriform algebras} and \ref{Example : Zinbiel algebras}. Then, in order to define such an algebraic structure on our cochain complex, we have to link shuffles and some sets of trunk maps. The link is provided by the map \eqref{definition of rhosigma}
$$
\rho : \mathrm{Sh}_{p_1,p_2} \to \mathrm{Hom}_{\mathbf{Set}}\big(N(X_{\mathbf{Tr}})(\square_{p_1+p_2}),N(X_{\mathbf{Tr}})(\square_{p_1}) \times N(X_{\mathbf{Tr}})(\square_{p_2})\big).
$$
which satisfies good properties with respect to the composition $\circ$ and the concatenation product $\star$ \eqref{concatenation product} defined on $\mathbb{S}_n$ (Lemma \ref{Lemma : graded dendriform algebra structure}). 
\paragraph{Section 5 : A graded associative algebra morphism from $H^{\bullet}(G,A)$ to $HR^{\bullet}(\mathrm{Conf}(G),A)$.} In this section we recall the definitions of group cohomology $H^{\bullet}(-)$ and of the cup product defined on it. Viewing a group $G$ as a category $G_{\mathbf{Cat}}$ and using the simplicial nerve $\mathrm{B}(G_{\mathbf{Cat}})$, we recall the simplicial definition of group cohomology (with trivial coefficient) and of its cup product. Now, on the one hand we have a cohomology theory based on simplicial relations and on the other hand a cohomology theory based on cubical relations. Therefore, in order to construct the expected morphism from $H^{\bullet}(G,A)$ to $HR^{\bullet}(\mathrm{Conf}(G),A)$, we have to find a relation between simplices and cubes. The construction of this morphism \big(Theorem \ref{Theorem : construction of the morphism S}\big) is based on the decomposition of the $n$-cubes into $n !$ $n$-simplex, that is, the existence of functors $\sigma : \Delta_n \to \square_n$ for all $\sigma \in \mathbb{S}_n$. We can represent geometrically such a functor in the following manner $(n=3,\sigma = (13))$ 
$$
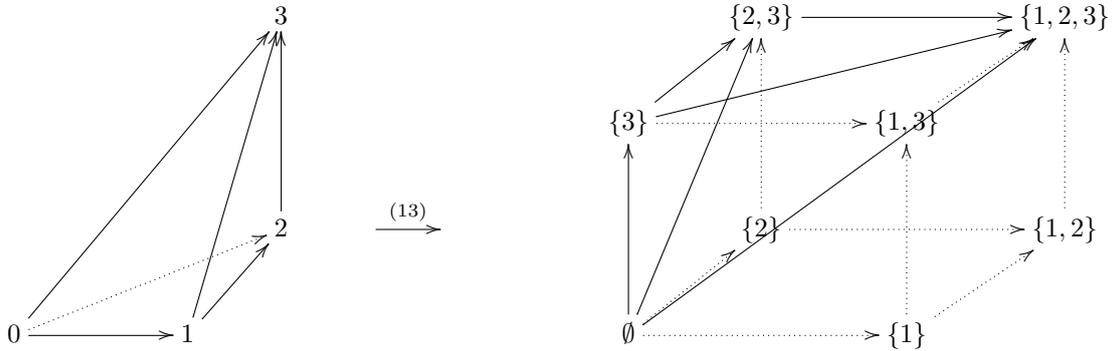
\xymatrix{
 & & & 3 &  & & & & \{2,3\} \ar[rr] & & \{1,2,3\}
\\
 & & & & & & & \{3\} \ar[ru] \ar@{.>}[rr]  \ar[urrr]& & \{1,3\} \ar@{.>}[ru] & 
\\
 & & & 2 \ar[uu] & \ar[r]^{(13)} & & & & \{2\} \ar@{.>}[uu] \ar@{.>}[rr] & & \{1,2\} \ar@{.>}[uu]
\\
0 \ar[rr] \ar[rrruuu] \ar@{.>}[rrru] &  & 1 \ar[ru] \ar[ruuu]  & & & & & \emptyset \ar[uu] \ar@{.>}[rr]  \ar@{.>}[ru] \ar[ruuu] \ar[uuurrr] & & \{1\}  \ar@{.>}[uu]  \ar@{.>}[ur] &
}
$$
We finish this section by proving that this morphism is compatible with the graded associative structures on $H^{\bullet}(-)$ and $HR^{\bullet}(-)$ \big(Theorem \ref{Theorem : a graded associative algebra morphism from group cohomology to rack cohomology}\big).
\section{Dendriform and Zinbiel algebras}
\paragraph{Shuffles.}
For $n \in \N$, we denote by $\mathbb{S}_n$ the group of permutations of $\{1,\dots,n\}$. For $p_1,p_2 \in \N$, a $(p_1,p_2)$-\textit{shuffle} is an element $\sigma$ of $\mathbb{S}_{p_1+p_2}$ satisfying 
\begin{align*}
\sigma(1) < \cdots < \sigma(p_1)  \, \text{ and } \, \sigma(p_1+1) < \cdots < \sigma(p_1+p_2).
\end{align*}
Remark that $\sigma$ satisfies $\sigma(p_1) = p_1+p_2$ or $\sigma(p_1+p_2) = p_1+p_2$. 
\par
Denote by $\Sh_{p_1,p_2}$ the set of $(p_1,p_2)$-shuffles, $\Sh_{p_1,p_2}^{p_1}$ the subset of elements $\sigma \in \Sh_{p_1,p_2}$ satisfying $\sigma(p_1)=p_1+p_2$ and $\Sh_{p_1,p_2}^{p_1+p_2}$ the subset of elements $\sigma \in \Sh_{p_1,p_2}$ satisfying $\sigma(p_1+p_2)= p_1+p_2$.
\par
Identifying a permutation $\sigma \in \S_n$ with the $n$-tuple $(\sigma(1),\dots,\sigma(n))$, the set $\Sh_{p_1,p_2}$ corresponds to the set of $n$-tuples $(a_1,\dots,a_{p_1+p_2})$ satisfying $a_1 < \dots < a_{p_1}$ and $a_{p_1+1} < \dots < a_{p_1+p_2}$. The cardinal of $\Sh_{p_1,p_2}$ (resp. $\Sh_{p_1,p_2}^{p_1}, \, \Sh_{p_1,p_2}^{p_1+p_2}$) is the binomial coefficient $\begin{pmatrix} p_1+p_2 \\ p_1 \end{pmatrix}$ (resp. $\begin{pmatrix} p_1+p_2 -1 \\ p_1-1 \end{pmatrix}, \, \begin{pmatrix} p_1+p_2-1 \\ p_1 \end{pmatrix}$).
\par
In the same way, given $p_1,p_2,p_3 \in \N$, the set of $(p_1,p_2,p_3)$-\textit{shuffles} is defined as the subset of elements $\sigma \in \S_{p_1+p_2+p_3}$ satisfying
\begin{align*}
\sigma(1) < \dots < \sigma(p_1)\, , \, \sigma(p_1+1) < \dots < \sigma(p_1+p_2) \, \text{ and } \, \sigma(p_1+p_2+1) < \dots < \sigma(p_1+p_2+p_3).
\end{align*}
\par
For $p_1,p_2,p_3 \in \N$, there are bijections
\begin{align}
\Sh_{p_1,p_2+p_3} \times \Sh_{p_2,p_3} \stackrel{\alpha}{\longrightarrow} \Sh_{p_1,p_2,p_3} \, \text{ and }  \Sh_{p_1+p_2,p_3} \times \Sh_{p_1,p_2} \stackrel{\beta}{\longrightarrow} \Sh_{p_1,p_2,p_3} \label{bijections between shuffles}
\end{align}
given by $\alpha(\sigma,\gamma) = \sigma \circ (1_{p_1} \star \gamma)$ and $\beta(\sigma,\gamma) = \sigma \circ (\gamma \star 1_{p_3})$ where $\star : \mathbb{S}_p \times \mathbb{S}_q \to \mathbb{S}_{p+q}$ is the map defined by
\begin{align} \label{concatenation product}
(\sigma \star \gamma)(k) := \left\{\begin{array}{ll}
						\sigma(k) & \textrm{if } 1 \leq k \leq p, \\
						p + \gamma(k-p) & \textrm{if } p+1 \leq k \leq p+q.
						\end{array}\right. 
\end{align}
In subset notation, the formulas for $\alpha$ and $\beta$ are
\begin{align*}
\alpha\big((a_1,\dots,a_{p_1+p_2+p_3}),(b_1,\dots,b_{p_2+p_3})\big) = (a_1,\dots,a_{p_1},a_{p_1+b_1},\dots,a_{p_1+b_{p_2+p_3}}) 
\end{align*}
and
\begin{align*}
\beta\big((a_1,\dots,a_{p_1+p_2+p_3}),(b_1,\dots,b_{p_1+p_2})\big) = (a_{b_1},\dots,a_{b_{p_1+p_2}},a_{p_1+p_2+1},\dots,a_{p_1+p_2+p_3}).
\end{align*}
As we will see in the examples \ref{Example : Dendriform algebras} and \ref{Example : Zinbiel algebras}, the notion of a shuffle is closely related to the notion of a dendriform algebra.
\paragraph{Dendriform algebras.} A \textnormal{(graded) dendriform algebra} is a (graded) vector space $D$ provided with two (graded) linear maps, 
$$
\succ,\prec : D \otimes D \to D,
$$
which satisfy for all $x,y,z \in D$ the following relations
\begin{align*}
x \succ (y \succ z) &= (x \succ y) \succ z + (x \prec y) \succ z \\
(x \succ y) \prec z &= x \succ (y \prec z) \\
(x \prec y) \prec z &= x \prec (y \prec z) + x \prec (y \succ z)
\end{align*}
\begin{example}{$T(V)$ and the free dendriform algebra $F(V)$}\label{Example : Dendriform algebras}
\begin{itemize}
\item Let $V$ be a vector space. We define a dendriform algebra structure on $T(V)$ by setting
$$
v_1\dots v_p \succ v_{p+1} \dots v_{p+q} := \displaystyle \sum_{\sigma \in Sh_{p,q}^{p+q}} v_{\sigma^{-1}(1)}\dots v_{\sigma^{-1}(p+q-1)}v_{p+q},
$$
and
$$
v_1\dots v_p \prec v_{p+1} \dots v_{p+q} := \displaystyle \sum_{\sigma \in Sh_{p,q}^{p}} v_{\sigma^{-1}(1)} \dots v_{\sigma^{-1}(p+q-1)}v_{p}.
$$
\item Let $V$ be a vector space. J.-L. Loday has defined a structure of dendriform algebra on $F(V) := \displaystyle \bigoplus_{n \geq 0} \K[Y_n] \otimes V^{\otimes n}$ where $Y_n$ is the set of planar binary trees with $n+1$ leaves (cf. \cite{LodayDialgebras}). This algebra is the free dendriform algebra associated to $V$.
\item The next example of a graded dendriform algebra uses the well known \textit{Pascal formula}.
\par
Let $A$ be an associative algebra and consider the graded $\K$-module $\displaystyle \bigoplus_{n \geq 0} A_n$ where $A_n = A$ for all $n \in \N$. We define the structure of a graded dendriform algebra on it by setting for $a \in A_p$ and $b \in A_q$
$$
a \succ b := \begin{pmatrix} p+q-1 \\ p \end{pmatrix} ab,
$$
and 
$$
a \prec b :=  \begin{pmatrix} p+q-1 \\ p-1 \end{pmatrix} ab.
$$
\end{itemize}
\end{example}
\paragraph{Relation with associative algebras.} A (graded) dendriform algebra is a particular example of a (graded) associative algebra. Indeed, let $(D,\succ,\prec)$ be a (graded) dendriform algebra, we define a new product $*$ on $D$ by setting
$$
x * y := x \succ y + x \prec y.
$$
\begin{proposition}
The product $*$ is associative.
\end{proposition}
Hence there exists a functor from the category of dendriform algebras to the category of associative algebras
$$
\mathbf{Dend} \to \mathbf{As}
$$
\par
The product $*$ is not necessarily commutative, but it becomes commutative if $x \succ y = y \prec x$. This condition leads us to the notion of a \textit{commutative dendriform algebra}, which is also called \textit{Zinbiel algebra}. 
\paragraph{Zinbiel algebras.}\label{Zinbiel algebra} A \textnormal{Zinbiel algebra} is a vector space $D$ provided with a linear map 
$$
\succ : D \otimes D \to D,
$$
which satisfies for all $x,y,z \in D$ the following relation
$$
x \succ (y \succ z) = (x \succ y) \succ z + (y \succ x) \succ z.
$$
Remark that the variables do not stay in the same order, thus in the graded case we have to be careful. A \textit{graded Zinbiel algebra} is a graded vector space $D$ provided with a graded linear map of degree $0$ $\succ :D \otimes D \to D$ which satisfies for all $x,y,z \in D$ the relation 
$$
x \succ (y \succ z) = (x \succ y) \succ z + (-1)^{pq} (y \succ x) \succ z,
$$
where $p$ is the degree of $x$ and $q$ is the degree of $y$.
\begin{example}{$T(V)$, $HL^{\bullet}(\mathfrak{g},A)$}\label{Example : Zinbiel algebras}
\begin{itemize}
\item Let $V$ be a vector space. There is a Zinbiel algebra structure defined on $T(V)$ by setting
$$
v_1\dots v_p \succ v_{p+1} \dots v_{p+q} := \displaystyle \sum_{\sigma \in Sh_{p,q}^{p+q}} v_{\sigma^{-1}(1)}\dots v_{\sigma^{-1}(p+q-1)}v_{p+q}.
$$
\item Let $\gg$ be a Leibniz algebra and $A$ a commutative algebra. In \cite{Lodaycup} J.-L. Loday defines a graded Zinbiel product on $HL^*(\gg,A)$, the Leibniz cohomology of $\gg$ with values in the trivial $\mathfrak{g}$-representation $A$.
\end{itemize}
\end{example}
\paragraph{Relation with dendriform algebras and commutative algebras.} As stated above, a Zinbiel algebra is kind of a commutative dendriform algebra. Indeed, let $(D,\succ)$ be a Zinbiel algebra, and define a second product $\prec$ by the formula $x \prec y := y \succ x$ for all $x,y \in D$.
\begin{proposition}
$(D,\succ,\prec)$ is a dendriform algebra, and $(D,*)$ is commutative. Conversely, a dendriform algebra $(D,\succ,\prec)$ where $x \succ y = y \prec x$ for all $x,y \in D$,  is a Zinbiel algebra.
\end{proposition}
Hence, there exist functors from the category of Zinbiel algebras to the category of commutative algebras, and to the category of dendriform algebras. All these functors fit into the following diagram
$$
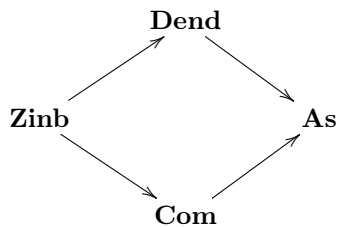
\xymatrix{ & \textbf{Dend} \ar[rd] & \\ \textbf{Zinb} \ar[ru] \ar[rd] & & \textbf{As} \\ & \textbf{Com} \ar[ru] & }
$$
\begin{remark}\label{Remark : beautiful diagram}
This diagram fits into a (beautiful) butterfly diagram (cf. \cite{LodayDialgebras})
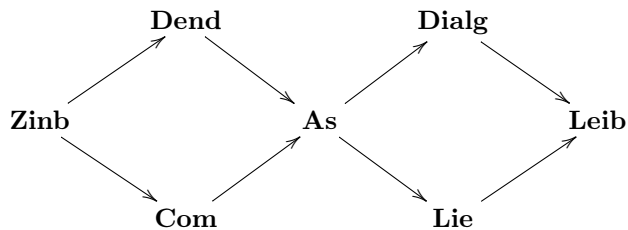
$$
\xymatrix{ & \textbf{Dend} \ar[rd] & & \textbf{Dialg} \ar[rd] & \\ \textbf{Zinb} \ar[ru] \ar[rd] & & \textbf{As} \ar[ru] \ar[rd] & & \textbf{Leib} \\ & \textbf{Com} \ar[ru] & & \textbf{Lie} \ar[ru] & }
$$
where \textbf{Dialg} is the category of dialgebras, \textbf{Leib} is the category of Leibniz algebras, and \textbf{Lie} is the category of Lie algebras. We can remark that each operad on the left (\textbf{Zinb, Dend, Com, As}) is Koszul dual to the operad on the right (\textbf{Leib, Dialg, Lie, As}) which is symmetric to it with respect to the vertical line passing through \textbf{As}.  
\end{remark}
\section{Racks}
Racks are sets provided with a product which encapsulates some properties of the conjugation in a group (especially the self-distributivity). This algebraic structure has been introduced in knot theory to define invariant for knots and links (cf. \cite{CarterSaito} or \cite{FennRourke}). We recall some basic definitions about racks and its cohomology theory, and the relations between the category of racks and the category of groups.
\paragraph{Shelves.} A \textit{shelf} is a set $X$ provided with a binary product $\rhd : X \times X \to X$ satisfying $x \rhd (y \rhd z) = (x \rhd y) \rhd (x \rhd z)$ for all $x,y,z \in X$. A \textit{pointed shelf} is the data of a shelf $(X,\rhd)$ provided with an element $1 \in X$, called the unit, satisfying $1 \rhd x = x$ and $x \rhd 1 = 1$ for all $x \in X$.
\vskip 0.2 cm
\noindent
\textbf{Notation :} Let $(X,\rhd)$ be a shelf. For all $x \in X$ denote by $c_x$ the map from $X$ to $X$ defined by $c_x(y) = x \rhd y$. Because in a shelf the product $\rhd$ is non associative, we have to be careful with the parenthesis. In the sequel the expression $(c_{x_1} \circ \dots \circ c_{x_{n_1}})(x_{n})$ (bracketing from the right) is denoted by $x_1 \conj \dots \conj x_n$ or $\displaystyle \prod_{1 \leq i \leq n} x_i$.
\vskip 0.2 cm
Let $X$ and $X'$ be two shelves. A \textit{morphism of shelves} from $X$ to $X'$ is a map $f : X \to X'$ which preserves products. If the shelves are pointed then a \textit{morphism of pointed shelves} is a morphism of shelves which preserves units.
\vskip 0.2 cm
The category of shelves (resp. pointed shelves) is denoted by $\mathbf{Shelf}$ (resp. $\mathbf{pShelf}$).
\paragraph{Racks.} A \textit{(pointed) rack} is a (pointed) shelf $(X,\rhd)$ where for all $x \in X$ the map $c_x : X \to X$ is a bijection. Let $X$ and $X'$ be two (pointed) racks. A \textnormal{morphism of (pointed) racks} from $X$ to $X'$ is a morphism of (pointed) shelves $f : X \to X'$.
\vskip 0.2 cm
The category of racks (resp. pointed racks) is denoted by $\mathbf{Rack}$ (resp. $\mathbf{pRack}$).
\begin{example}\label{Example : groups and augmented racks}{Groups and augmented racks}
\begin{itemize}
\item Let $G$ be a group. Define a rack structure on the set $G$ by $x \conj y = xyx^{-1}$. This rack is denoted $Conj(G)$. It is pointed by $e$, the unit of $G$.  
\item Let $G$ be a group, $X$ be a $G$-set and $f : X \to G$ a $G$-map (where $G$ acts on itself by conjugation). Define a rack structure on $X$ by $x \conj y = f(x) \cdot y$. If there exists a fixed point $1 \in X$ such that $f(1) = e$, then this rack is pointed by $1$. 
\end{itemize}
\end{example}
Example \ref{Example : groups and augmented racks} defines a functor $\mathrm{Conj} : \mathbf{Grp} \to \mathbf{Rack}$ from the category of groups to the category of racks. There is a functor $\mathrm{As} :  \mathbf{Rack} \to \mathbf{Grp}$ defined by
\begin{align*}
\mathrm{As}(X) &:= \mathcal{F}(X) / < x \conj y = xyx^{-1} >,\\
\mathrm{As}(f) &:= \overline{\mathcal{F}(f)}, 
\end{align*}
where $\mathcal{F}(X)$ is the free group generated by the set $X$, and $\overline{\mathcal{F}(f)}$ is the map induced by $\mathcal{F}(f)$ by passing to quotients.
\begin{proposition} \label{As is left adjoint to Conj}
The functor $\mathrm{As}$ is left adjoint to the functor $\mathrm{Conj}$. 
\end{proposition}
\paragraph{Rack modules.} Let $X$ be a rack. A \textit{left $X$-module} is an abelian group $A$ provided with a right linear map 
\begin{align*}
\cdot : X \times A \to A; (x,a) \mapsto x \cdot a
\end{align*}
satisfying  $x \cdot ( y \cdot a) = (x \conj y) \cdot (x \cdot a)$ for all $a \in A$ and $x,y \in X$. If the rack is pointed then we demand the additional axiom : $1 \cdot a = a, \, \forall a \in A$. Let $A$ and $B$ be two left $X$-modules. A map $f : A \to B$ is a morphism of left $X$-modules if $f(x \cdot a) = x \cdot f(a)$ for all $(a,x) \in A \times X$. We denote by $X-\textbf{Mod}$ the category of left $X$-modules.
\begin{example}\label{module over a group is a module over a rack}
Let $G$ be a group and $A$ be a left $G$-module, then $A$ is a left $\mathrm{Conj}(G)$-module.
\end{example}
\begin{remark} \label{equivalence between X-Mod and As(X)-Mod}
Let $X$ be a rack and $A$ be a left $X$-module. There is an equivalence of categories between $X-\mathbf{Mod}$ and $\mathrm{As}(X)-\mathbf{Mod}$.
\end{remark}
\paragraph{Rack cohomology.} \label{Definition : rack cohomology} Let $X$ be a rack and $A$ be a left $X$-module. We define a cochain complex $\{CR^n(X,A),d_R^{n}\}_{n \in \mathbb{N}} $ by
\begin{align*}\displaystyle
CR^n(X,A) &:= \mathrm{Hom}_{\mathbf{Set}}(X^n,A),\\
d_R^{n+1} f &:= \sum_{i=1}^{n+1} (-1)^i \big(d_{i,0}^{n+1} f - d_{i,1}^{n+1} f \big),
\end{align*}
where
\begin{align*}
d_{i,\epsilon}^{n+1} f (x_1,\dots,x_{n+1}):= \left\{ \begin{array}{ll}
					f(x_1,\dots,x_{i-1},x_i \conj x_{i+1},\dots,x_i \conj x_{n+1}) & \textrm{ if } \epsilon = 0,\\
					(x_1 \conj \dots \conj x_i) \cdot f(x_1,\dots,x_{i-1},x_{i+1},\dots,x_{n+1}) & \textrm{ if } \epsilon = 1.
					\end{array} \right.
\end{align*}
\par
The fact that $d_R$ is a differential comes from the \textit{cubical identities} satisfied by the family of maps $\{d_{i,\epsilon}^{n}\}$, that is, for all $1 \leq i < j \leq n+1$ and $\epsilon,\omega \in \{0,1\}$, we have the identities :
$$
d_{i,\epsilon}^{n+1} \circ d_{j-1,\omega}^n = d_{j,\omega}^{n+1} \circ d_{i,\epsilon}^n.
$$
The cohomology associated to this cochain complex is called the \textit{rack cohomology of $X$ with coefficients in $A$} and denoted $HR^{\bullet}(X,A)$.
\begin{remark} 
In this paper we will consider the rack cohomology with trivial coefficients, that is when $A$ is a trivial left $X$-module. In this case we have $d_{n+1,0}^{n+1} f = d_{n+1,1}^{n+1} f$ for all $f \in CR^{n}(X,A)$. 
\end{remark}
\begin{remark}
If $X$ is a pointed rack, then a subcomplex $CR_p^{\bullet}(X,A) = \{CR^n_p(X,A),d_R^n\}_{n \geq 1}$ is defined by
$$
CR^n_p(X,A) := \{f \in \mathrm{Hom}_{\mathbf{Set}}(X^{n},A) \, | \, f(x_1,\dots,1,\dots,x_n) = 0\}.
$$
The cohomology associated to this cochain complex is called the \textit{pointed rack cohomology of $X$ with coefficients in $A$} and denoted $HR^{\bullet}_p(X,A)$.
\end{remark}
\begin{remark}
In the definition of $\{CR^n(X,A),d_R^n\}_{n \in \mathbb{N}}$ (resp. $\{CR_p^n(X,A),d_R^n\}_{n \in \mathbb{N}}$ we don't use the second axiom of a rack, that is the bijectivity of $c_x$ for all $x \in X$. Hence this cochain complex is well defined for a shelf (resp. pointed shelf). It correponds to the cochain complex associated to the \textit{multi-shelf} $(X,\rhd,\lhd)$, where $x \lhd y = x$, defined in \cite{Przytycki_Sikora_Distributive_Products_And_Their_Homology}. 
\end{remark}
\section{The category of trunks}
Trunks are the fundamental tools in order to solve our problems combinatorially. This section is based on \cite{FennRourkeSanderson}. Note that our definitions are slightly different from theirs for we are working in the left rack context instead of the right rack context. The two theories are obviously equivalent. 
\subsection{Definitions and examples}
\paragraph{Categories.} Following \cite{Maclane_Categories_for_the_working_mathematician}, let us recall the definition of a (small) category as a graph provided with extra properties. Let $\mathcal{G} = (V,E,s,t)$ be a directed graph with $V$ the set of vertices, $E$ the set of edges, $s :E \to V$ the source map and $t : E \to V$ the target map. We define a subset of $E \times E \times E$ by
$$
E \times_V E \times_V E := \{(a,b,c) \in E \times E \times E \, | \, t(a) = s(b), \, t(b) = s(c)\},
$$
and subsets of $E \times E$ by 
\begin{align*}
E \times_V E &:= \{(a,b) \in E \times E \, | \, t(a) = s(b)\}, \\
V \times_V E &:= \{(s(a),a) \in V \times E \, | \, a \in E\}, \\ 
E \times_V V &:= \{(a,t(a)) \in E \times V \, | \, a \in E\}.
\end{align*}
A \textit{category} is a directed graph $\mathcal{C} = (V,E,s,t)$ provided with a \textit{composition map} $c:E \times_V E \to E$ and an \textit{identity map} $i:V \to E$ satisfying the following commutative diagrams.
$$
\xymatrix{
E \times_V E \times_V E \ar[d]_{1 \times c} \ar[r]^{\quad c \times 1} & E \times_V E \ar[d]^{c} & V \times_V E \ar[r]^{i \times 1} \ar[rd]^{\pi_2} & E \times_V E \ar[d]^{c} & E \times_V V \ar[l]_{1 \times i} \ar[ld]_{\pi_1} \\ 
E \times_V E \ar[r]^{c} & E & & E & \\
}
$$
\par
The left diagram states the \textit{associativity} of the composition $c$ and the right diagram states that each $i(x)$ is a \textit{left and right unit} for the composition. 
\par
\begin{example}{Simplexes, cubes and groups.}\label{Example : category structure on square_n}
\begin{itemize}
\item Let $n \in \N$, and let us consider the directed graph $\Delta_n = \{V(\Delta_n),E(\Delta_n),s,t\}$ where $V(\Delta_n) := \{0,\dots,n\}, \, E(\Delta_n) := \{(i,j) \in V(\Delta_n) \times V(\Delta_n) \, | \, i \leq j\}, \, s(i,j) = i$ and $t(i,j) = j$. Define a composition map $c$ by $c\big((i,j),(j,k)\big) = (i,k)$ and an identity map $i$ by $i(j) = (j,j)$. In this way, the standard simplex $\Delta_n$ gives rise to a category denoted $\Delta_n$.
\item Let $n \in \N$, let us consider the directed graph $\square_n = \{V(\square_n),E(\square_n),s,t\}$ where $V(\square_n) := \{A \subseteq \{1,\dots,n\}\}, \, E(\square_n) := \{(A,B) \in V(\square_n) \times V(\square_n) \, | \, A \subseteq B\}, s(A,B) = A$ and $t(A,B) = B$. Define a composition map $c$ by $c\big((A,B),(B,C)\big) = (A,C)$ and an identity map $i$ by $i(A) = (A,A)$. In this way, the standard cube $\square_n$ gives rise to a category denoted $\square_n$.
\item Let $G$ be a group. Consider the directed graph $G = (\{\star\},G,s,t)$ with $s = t: G \to \{\star\}$, the canonical map from $G$ to the terminal object $\{\star\}$ in $\mathbf{Set}$. Define a composition map $c : G \times G \to G$ by $c(a,b) = \mu(a,b)$ (the product of $a$ and $b$ in $G$) and an identity map $i : \{\star\} \to G$ by $i(\star) = e$ (the unit in $G$). In this way, a group $G$ give rise to a category.
\end{itemize}
\end{example}
\paragraph{Functors.} Let $\mathcal{C}$ and $\mathcal{C}'$ be two categories. A \textit{functor} from $\mathcal{C}$ to $\mathcal{C}'$ is a graph map $F : \mathcal{C} \to \mathcal{C}'$ which preserves composition and units, i.e $F$ satisfies the following commutative diagrams
$$
\xymatrix{ E \times_V E \ar[r]^{c} \ar[d]_{F \times F} & E \ar[d]^{F} & & V \ar[r]^{i} \ar[d]_{F} & E \ar[d]^{F} 
\\
E' \times_V E'  \ar[r]^{c'} & E' & & V' \ar[r]^{i'} & E'}
$$
\begin{example}{Faces of simplexes, faces of cubes and permutations.}\label{Example : functors}
\begin{itemize}
\item Let $n \in \mathbb{N}$. For all $0 \leq i \leq n$ define a functor $\partial_i^n$ from the category $\Delta_{n-1}$ to the category$\Delta_n$ by
$$
\partial_i^n(k) := \left\{\begin{array}{ll}
k & \textrm{if } 0 \leq k \leq i-1,\\
k+1 & \textrm{if } i \leq k \leq n-1.
\end{array}
\right.
$$ 
Given two vertices in $\Delta_{n-1}$ there is exactly one edge between them, hence there is only one way to define $\partial_i^n$ on the set of edges in order to have a graph map.
\item Let $n \in \mathbb{N}$. For all $(i,\epsilon)\in \{1,\dots,n\} \times \{0,1\}$ define a functor $\partial_{i,\epsilon}^n$ from the category $\square_{n-1}$ to the category $\square_n$ by
$$
\partial_{i,\epsilon}^n(A) = \left\{ \begin{array}{ll} 
A_{ < i} \amalg t_{+1}(A_{\geq i}) & \textrm{ if } \epsilon = 0,\\
A_{ < i} \amalg t_{+1}(A_{\geq i}) \amalg \{i\} & \textrm{ if } \epsilon = 1.
\end{array}\right. 
$$
where $A_{< i} = \{a \in A \, | \, a < i\}, \, A_{\geq i} = \{a \in A \, | \, a \geq i\}$ and $t_{+1}(B) = \{b+1 \, | \, b \in B\}$. Given two vertices in $\square_{n-1}$ there is at most one edges between them, hence there is only one way to define $\partial_{i,\epsilon}^n$ on the set of edges in order to have a graph map.
\item Let $n \in \mathbb{N}$. For all $\sigma \in \mathbb{S}_n$ define a functor $\sigma$ from the category $\Delta_n$ to the category $\square_n$ by
$$
\sigma(0) := \emptyset \, \text{ and } \, \sigma(k) := \{\sigma(1),\dots,\sigma(k)\} \, \forall k \geq 1.
$$
\end{itemize}
\end{example}
As explained in the introduction, in order to define trunks we want to replace the associativity relation for the composition in a category by the \textit{self-distributivity} relation.
\paragraph{Trunks.}
Let $\mathcal{G} = (V,E,s,t)$ be a directed graph. We denote by $\mathrm{S}(\mathcal{G})$ the subset of $E \times_V E \times E \times_V E$ of elements $(a,b,c,d)$ satisfying $s(a) = s(c), \, t(a) = s(b), \, t(c) = s(d)$ and $t(b) = t(d)$. A \textit{trunk} is a directed graph $\mathcal{T} = (V,E,s,t)$ provided with a subset $\Gamma$ of $\mathrm{S}(\mathcal{T})$. An element in $\Gamma$ is called a \textit{preferred square}.
\par
A \textit{pointed trunk} is a trunk $(\mathcal{T}=(V,E,s,t),\Gamma)$ together with an identity map $i:V \to E$ satisfying $\big(a,i(t(a)),i(s(a)),a\big) \in \Gamma$ and $\big(i(s(a)),a,a,i(t(a))\big) \in \Gamma$ for all $a \in E$. 
\begin{example}{Categories and cubes.}\label{Example : categories and cubes are trunks}
\begin{itemize}
\item Let $\mathcal{C} = (V,E,s,t)$ be a category. We define a trunk $\mathrm{Tr}(\mathcal{C})$ by the pair $(\mathcal{C},\Gamma)$ with $\Gamma$ the set of commutative diagrams in $\mathcal{C}$ 
$$
\Gamma = \{(a,b,c,d) \in \mathrm{S}(E) \, | \, c(a,b) = c(c,d)\}.  
$$
The identity map $i$ in the category $\mathcal{C}$ provides $\mathrm{Tr}(\mathcal{C})$ with a pointed trunk structure.
\item Let $n \in \mathbb{N}$. Consider the graph $\square_n$ with set of vertices $V(\square_n) = \{A \subseteq \{1,\dots,n\}\}$ and set of edges $E(\square_n) = \{(A,A \amalg \{k\}) \, | \, k \notin A  \}$. Take as set of preferred square $\Gamma_{\square_n}$ the set equals to 
$$
\{\big((A,A \amalg \{k\}),(A \amalg \{k\},A \amalg \{k,l\}),(A,A \amalg \{l\}),(A \amalg \{l\},A \amalg\{k,l\})\big) \, | \,  k < l\}.
$$
Then, in this way, the standard cube $\square_n$ give rise to a trunk denoted $\square_n$.
\end{itemize}
\end{example}
\paragraph{Corner trunks.}
A \textit{corner trunk} is a trunk $(\mathcal{T},\Gamma)$ where $\Gamma$ is the graph of a map $c : E \times_V E \to E \times_V E$ satisfying the following commutative diagram
$$
\xymatrix{ E \times_V E \times_V E \ar[d]_{c \times 1} & E \times_V E \times_V E \ar[l]_{1 \times c} \ar[r]^{c \times 1} & E \times_V E \times_V E \ar[d]^{1 \times c} \\
E \times_V E \times_V E \ar[r]^{1 \times c} & E \times_V E \times_V E  & E \times_V E \times_V E \ar[l]_{c \times 1}}
$$
The map $c$ is called the \textit{composition map}. The relation satisfied by $c$ can be decomposed into three relations, called \textit{bidistributivity relations}, according to the projections onto the three components in $E \times_V \times_V E$. To describe these relations, let us denote by $\rhd$ (resp. $\lhd$) the composition $\mathrm{pr}_1 \circ c$ (resp. $\mathrm{pr}_2 \circ c)$. This defines maps $\rhd$ and $\lhd$ from $E \times_V E$ to $E$, and the relations are :
\begin{itemize}
\item[$\bullet$] $a \rhd (b \rhd c) = (a \rhd b) \rhd ((a \lhd b) \rhd c)$, 
\item[$\bullet$] $(a \rhd b) \lhd ((a \lhd b) \rhd c) = (a \lhd (b \rhd c)) \rhd (b \lhd c)$,
\item[$\bullet$] $(a \lhd b) \lhd c = (a \lhd (b \rhd c)) \lhd (b \lhd c)$,
\end{itemize} 
for all $(a,b,c) \in E \times_V E \times_V E$.
\par
In case the trunk is pointed, we have necessarily the following relations, called \textit{unit relations} :
\begin{itemize}
\item $a \rhd i(t(a)) = i(s(a)) \, \text{ and } \, i(s(a) \rhd a = a$,
\item $a \lhd i(t(a)) = a \, \text{ and } \, i(s(a)) \lhd a = i(t(a))$.
\end{itemize}
\begin{example}{Racks and $\mathbb{N}$.}\label{Example : racks and natural numbers are trunks}
\begin{itemize}
\item Let $(X,\rhd)$ be a rack or a shelf. Consider the directed graph $\mathcal{T}_X = (\{\star\},X,s,t)$ with $s=t : X \to \{\star\}$, the unique map from $X$ to the terminal object in $\mathbf{Set}$. Define a composition map $c : X^2 \to X^2$ by $c(a,b) = (a \rhd b,a)$. The set of preferred square $\Gamma_X$ is equal to $\mathrm{S}(X) = X^4$. In this case, because the map $\lhd$ is the first projection, the second and third bidistributivity relations are trivials and the first is equivalent to the self-distributivity of $\rhd$. Remark that $\mathcal{T}_X$ is a pointed corner trunk if and only if $X$ is a pointed rack/shelf.
\item Consider the graph $\mathbb{N} = (\{\star\},\mathbb{N},s,t)$ with $s=t : \mathbb{N} \to \{\star\}$, the unique map from $\mathbb{N}$ to the terminal object in $\mathbf{Set}$. Define a composition map $c : \mathbb{N}^2 \to \mathbb{N}^2$ by $c(i,j) = (\max(i,j),\min(i,j))$. The set of preferred square $\Gamma_{\mathbb{N}}$ equals to $\mathrm{S}(\mathbb{N}) = \mathbb{N}^4$. Then, in this way, the set of natural numbers $\mathbb{N}$ give rise to a corner trunk.
\end{itemize}
\end{example}
\paragraph{Trunk maps.} Let $(\mathcal{T},\Gamma)$ and $(\mathcal{T}',\Gamma')$ be two trunks. A \textit{trunk map} from $\mathcal{T}$ to $\mathcal{T}'$ is a graph map $F : \mathcal{T} \to \mathcal{T}'$ mapping $\Gamma$ to $\Gamma'$. If the trunks are pointed then $F$ is a \textit{pointed trunk map} if $F$ preserves the units. In case the trunks $\mathcal{T}$ and $\mathcal{T'}$ are corner trunks, $F : \mathcal{T} \to \mathcal{T'}$ is a trunk map if and only if $F$ preserves the composition map $c$, or equivalently, if and only if $F$ preserves the products $\rhd$ and $\lhd$.
\begin{example}{Functors, faces of a cube and shuffles.}\label{Example : functors, faces and degeneracies are trunk maps} 
\begin{itemize}
\item Let $F$ be a functor from a category $\mathcal{C}$ to a category $\mathcal{C}'$. Because $F$ maps a commutative diagram in $\mathcal{C}$ to a commutative diagram in $\mathcal{C}'$, it induces a pointed trunk map $\mathrm{Tr}(F)$ from the pointed trunk $\mathrm{Tr}(\mathcal{C})$ to $\mathrm{Tr}(\mathcal{C})$ (cf. Example \ref{Example : categories and cubes are trunks}). 
\item Let $n \in \mathbb{N}$ and consider the graph map $\partial_{i,\epsilon}^n : \square_{n-1} \to \square_n$ defined in Example \ref{Example : functors}. 
\par
Given $1\leq k < l \leq n-1$, we have $\partial_{i,0}^{n}(\{k\}) < \partial_{i,0}^n(\{l\})$. Moreover for all $A, B$ such that $A \cap B = \emptyset$ we have $\partial_{i,\epsilon}^n(A \amalg B) = \partial_{i,\epsilon}^n(A)  \amalg \partial_{i,0}^n(B)$. Thus $\partial_{i,\epsilon}^n$ maps a preferred square in $\square_{n-1}$ to a preferred square in $\square_n$. Therefore, the maps $\partial_{i,\epsilon}^n$ are trunk maps.
\item Let $p_1,p_2 \in \mathbb{N}$ and $\sigma \in \mathrm{Sh}_{p_1,p_2}$. The permutation $\sigma$ induces a graph map $\sigma : \square_{p_1+p_2} \to \square_{p_1+p_2}$ on the vertices by 
$$
\sigma(\{a_1,\dots,a_k\}) := \{\sigma(a_1),\dots,\sigma(a_k)\}.
$$
Define two graph maps $\sigma \circ i_{p_1} : \square_{p_1} \to \square_{p_1+p_2}$ and $\sigma \circ i_{p_2} : \square_{p_2} \to \square_{p_1+p_2}$ where
\begin{align*}
i_{p_1} &:= \partial_{p_1+p_2,0}^{p_1+p_2} \circ \dots \circ \partial_{p_1+1,0}^{p_1+1},\\
i_{p_2} &:= \partial_{1,1}^{p_1+p_2} \circ \dots \circ \partial_{1,1}^{p_2+1}.
\end{align*}
We have
\begin{align*}
(\sigma \circ i_{p_1})(A) &= \{\sigma(a_1),\dots,\sigma(a_k)\},\\
(\sigma \circ i_{p_2})(A) &= \{\sigma(1),\dots,\sigma(p_1),\sigma(p_1+a_1),\dots,\sigma(p_1+a_k)\}. 
\end{align*}
Let $1 \leq k < l \leq p_1$, $\sigma \in \mathrm{Sh}_{p_1,p_2}$ implies $\sigma(k) < \sigma(l)$. Moreover for all $A,B$ such that $A \cap B = \emptyset$ we have $(\sigma \circ i_{p_1})(A \amalg B) = (\sigma \circ i_{p_1})(A) \amalg (\sigma \circ i_{p_1})(B)$, thus $(\sigma \circ i_{p_1})(A \amalg \{k\}) = (\sigma \circ i_{p_1})(A) \amalg \{\sigma(k)\}$ and $(\sigma \circ i_{p_1})(A \amalg \{l\}) = (\sigma \circ i_{p_1})(A) \amalg \{\sigma(l)\}$ with $\sigma(k) < \sigma(l)$. Thus $\sigma \circ i_{p_1}$ maps a preferred square in $\square_{p_1}$ to a preferred square in $\square_{p_1+p_2}$ and $\sigma \circ i_{p_1} : \square_{p_1} \to \square_{p_1+p_2}$ is a trunk map.
\par
In the same way, let $1 \leq k < l \leq p_2$, $\sigma \in \mathrm{Sh}_{p_1,p_2}$ implies $\sigma(p_1+k) < \sigma(p_1+l)$. Moreover for all $A,B$ such that $A \cap B = \emptyset$ we have $(\sigma \circ i_{p_2})(A \amalg B) = (\sigma \circ i_{p_2})(A) \amalg (\sigma \circ i_{p_2})(B)$, thus $(\sigma \circ i_{p_2})(A \amalg \{k\}) = (\sigma \circ i_{p_2})(A) \amalg \{\sigma(p_1+k)\}$ and $(\sigma \circ i_{p_2})(A \amalg \{l\}) = (\sigma \circ i_{p_2})(A) \amalg \{\sigma(p_1+l)\}$ with $\sigma(p_1+k) < \sigma(p_1+l)$. Thus $\sigma_ \circ i_{p_2}$ maps a preferred square in $\square_{p_2}$ to a preferred square in $\square_{p_1+p_2}$ and $\sigma \circ i_{p_2} : \square_{p_2} \to \square_{p_1+p_2}$ is a trunk map.
\par
Using these maps, we will define the graded dendriform algebra product on rack/shelf cohomology.
\end{itemize}
\end{example} 
The category of (pointed) trunks is the category with the (pointed) trunks as vertices and the (pointed) trunk maps as edges. This category is denoted by $\mathbf{Trunk}$ (resp. $\mathbf{pTrunk}$).  
\paragraph{Relations between categories and trunks.} In Example \ref{Example : categories and cubes are trunks} we have seen that a category $\mathcal{C}$ gives rise to a trunk $\mathrm{Tr}(\mathcal{C})$. Moreover we have seen in Example \ref{Example : functors, faces and degeneracies are trunk maps} that a functor $F$ from a category $\mathcal{C}$ to a category $\mathcal{C}'$ determines a trunk map $\mathrm{Tr}(F)$ from $\mathrm{Tr}(\mathcal{C})$ to $\mathrm{Tr}(\mathcal{C}')$. Hence a functor $\mathrm{Tr} : \mathbf{Cat} \to \mathbf{Trunk}$, from the category of categories $\mathbf{Cat}$ to the category of trunks $\mathbf{Trunk}$ defined by
\begin{align*}
\mathbf{Cat} &\stackrel{\mathrm{Tr}}{\longrightarrow} \mathbf{Trunk}\\
\mathcal{C} &\longmapsto \mathrm{Tr}(\mathcal{C})\\
F & \longmapsto \mathrm{Tr}(F) 
\end{align*}
\par
Let $(\mathcal{T},\Gamma)$ be a trunk. Consider the free category $\mathcal{F}(\mathcal{T})$ generated by the graph $\mathcal{T}$. For all pairs of objects $(x,y)$ in $\mathcal{F}(\mathcal{T})$ we define a binary relation $R_{x,y}$ on $\mathrm{Hom}_{\mathcal{F}(\mathcal{T})}(x,y)$ by $(a_1,\dots,a_n) R_{x,y} (b_1,\dots,b_m)$ if and only the two following points are satisfied :
\begin{enumerate}
\item $n=m$,
\item For all $1 \leq k \leq n$ such that $a_k \neq b_k$ then $(a_k,a_{k+1},b_k,b_{k+1}) \in \Gamma$ or $(a_{k-1},a_{k},b_{k-1},b_k) \in \Gamma$.
\end{enumerate}
We define a category $\mathrm{Cat}(\mathcal{T})$ by taking the quotient of $\mathcal{F}(\mathcal{T})$ by the binary relation $R$ (cf. \cite{Maclane_Categories_for_the_working_mathematician} for the definition of a quotient category). Then the category $\mathrm{Cat}(\mathcal{T})$ has the same vertices as $\mathcal{T}$ and the set of edges is the set of strings of composable edges in $\mathcal{T}$ where two strings of edges are identified if they have the same length and are equal "up to preferred squares in $\Gamma$". Hence 
all commutative diagrams in $\mathrm{Cat}(\mathcal{T})$ comes from preferred squares in $(\mathcal{T},\Gamma)$.
$$
x_0 \stackrel{a_0}{\rightarrow} x_1 \stackrel{a_1}{\rightarrow} \dots \stackrel{a_n}{\rightarrow} x_{n+1} \sim x_0 \stackrel{b_0}{\rightarrow} x_1 \stackrel{b_1}{\rightarrow} \dots \stackrel{b_n}{\rightarrow} x_{n+1}  
$$
if and only if
$$
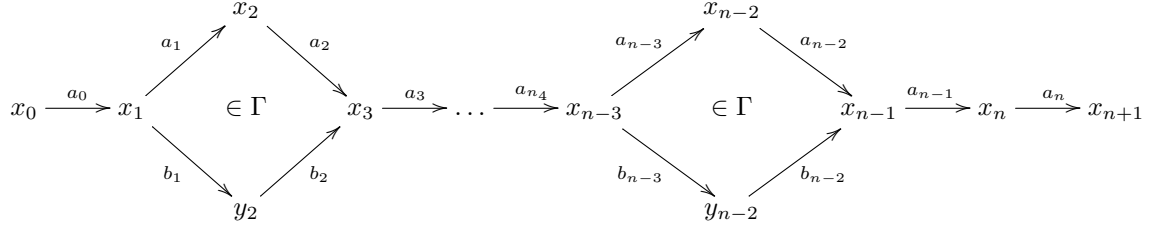
\xymatrix{
& & x_2 \ar[rd]^{a_2} & & & & x_{n-2} \ar[rd]^{a_{n-2}} & & &
\\
x_0 \ar[r]^{a_0} & x_1 \ar[ru]^{a_1} \ar[rd]_{b_1} & \in \Gamma & x_3 \ar[r]^{a_3} & \dots \ar[r]^{a_{n_4}}& x_{n-3} \ar[ru]^{a_{n-3}} \ar[rd]_{b_{n-3}} & \in \Gamma & x_{n-1} \ar[r]^{a_{n-1}} & x_n \ar[r]^{a_n} & x_{n+1}
\\
& & y_2 \ar[ru]_{b_2} & & & & y_{n-2} \ar[ru]_{b_{n-2}} & & &
}
$$
\par
Let $F$ be a trunk map from $(\mathcal{T},\Gamma)$ to $(\mathcal{T}',\Gamma')$. As $F$ is a graph map, $\mathcal{F}(F)$ is a functor from $\mathcal{F}(\mathcal{T})$ to $\mathcal{F}(\mathcal{T'})$. As $F$ maps preferred squares in $\Gamma$ to preferred squares in $\Gamma'$, $\mathrm{Cat}(F)$ is a well defined functor from $\mathrm{Cat}(\mathcal{T})$ to $\mathrm{Cat}(\mathcal{T'})$.
\par
Hence, there is a functor $\mathrm{Cat} : \mathbf{Trunk} \to \mathbf{Cat}$, from the category of trunks to the category of (small) categories, defined by
\begin{align*}
\mathbf{Trunk} &\stackrel{\mathrm{Cat}}{\longrightarrow} \mathbf{Cat}\\
\mathcal{T} &\longmapsto \mathrm{Cat}(\mathcal{T})\\
F & \longmapsto \mathrm{Cat}(F) 
\end{align*}
and $\mathrm{Cat}$ is left adjoint to $\mathrm{Trunk}$.
\begin{align}
\mathrm{Cat} \stackrel{\theta}{\vdash} \mathrm{Tr}. \label{adjunction Cat dashv Tr}
\end{align}
\begin{remark}
For all $n \in \mathbb{N}$ we can consider $\square_n$ as a category (Example \ref{Example : category structure on square_n}) or as a trunk (Example \ref{Example : categories and cubes are trunks}). We have 
\begin{align} \label{Cat(square)=square}
\mathrm{Cat}(\square_n) = \square_n,
\end{align}
where on the left $\square_n$ is seen as a trunk and on the right as a category, and 
$$
\mathrm{Tr}(\square_n) \hookrightarrow \square_n,
$$   
where on the left $\square_n$ is seen as a category and on the right as a trunk.
\end{remark}
\subsection{Nerve of a trunk}
\paragraph{The category $\square$.} Objects of the \textit{the cubical category} $\square$ are the graphs $\square_n$ defined in Example \ref{Example : category structure on square_n}. Given $n \in \mathbb{N}$ vertices of $\square_n$ has been defined as the subsets of $\{1,\dots,n\}$. The bijection between the set of subsets of $\{1,\dots,n\}$ and the set $\{0,1\}^n$, given by $A \mapsto (\epsilon_1,\dots,\epsilon_n)$ with $\epsilon_k = 1$ if $k \in A$ and $0$ if not, allows us to describe $\square_n$ in a second way. In the first case, we will say that $\square_n$ is described in subset notation, and in the second case, in coordinate notation.
\par
For all $m,n \in \mathbb{N}$, define $\mathrm{Hom}_{\square}(\square_m,\square_n)$ as the subset of $\mathrm{Hom}_{\mathbf{Graph}}(\square_m,\square_n)$ generated by the face maps $\partial_{i,\epsilon}^k$ defined in Example \ref{Example : functors, faces and degeneracies are trunk maps}. In other words a morphism $\mathrm{Hom}_{\square}(\square_m,\square_n)$ is a graph map from $\square_m$ to $\square_n$ which is a composition of face maps. Using the description of $\square_n$ in terms of coordinates, $\partial_{i,\epsilon}^n$ is the graph map defined by
\begin{align*}
\partial_{i,\epsilon}^n(\epsilon_1,\dots,\epsilon_{n-1}) &= (\epsilon_1,\dots,\epsilon_{i-1},\epsilon,\epsilon_{i},\dots,\epsilon_{n-1}).
\end{align*}
The next proposition states that face maps satisfy relations called \textit{cubical relations}. This allows us to write a morphism in $\mathrm{Hom}_{\square}(\square_m,\square_n)$ in a canonical way.
\begin{proposition}[Cubical relations] \label{Proposition : Cubical relations} For all $1 \leq i < j \leq k+1$, and $\epsilon,\omega \in \{0,1\}$ we have the following relation :
\begin{align}
\partial_{i,\epsilon}^{k+1} \circ \partial_{j-1,\omega}^k &= \partial_{j,\omega}^{k+1} \circ \partial_{i,\epsilon}^k.
\end{align}
\end{proposition}
\begin{corollary}
Let $m < n \in \mathbb{N}$. Each element $f$ in $\mathrm{Hom}_{\square}(\square_m,\square_n)$ can be written uniquely as 
$$
f = \partial_{i_{n-m},\epsilon_{n-m}}^{n} \circ \dots \circ \partial_{i_1,\epsilon_1}^{m+1},
$$
with $i_1 < i_2 < \cdots < i_{n-m}$.
\end{corollary}
\paragraph{Cubical object in a category.} Let $\mathcal{C}$ be a category. A \textit{precubical object in $\mathcal{C}$} is a functor $N : \square^{op} \to \mathcal{C}$. Dually a \textit{precocubical object in $\mathcal{C}$} is a functor $N : \square \to \mathcal{C}$.
\par
Let $\mathbb{K}$ be a commutative ring and $\mathcal{C}$ be the category $\mathbf{Mod}_{\mathbb{K}}$ of modules over $\mathbb{K}$. From a precocubical $\mathbb{K}$-module $M$ (that is a precocubical object in $\mathbf{Mod}_{\mathbb{K}}$), define a cochain complex $M^{\bullet} = \{M^n,d^n\}$ by   
$$
M^n := M(\square_n)\, , \, d^n := \displaystyle \sum_{i=1}^{n} (-1)^i (M(\partial_{i,0}^n) - M(\partial_{i,1}^n)),
$$
The graded map $d^n$ is a differential thanks to the cubical relations (Proposition \ref{Proposition : Cubical relations}).
\begin{remark} \label{Remark : from precubical set to precocubical module}
Let $A$ be a $\mathbb{K}$-module. From a precubical set $S : \square^{op} \to \mathbf{Set}$ we can always define a precocubical $\mathbb{K}$-module by postcomposed $S$ by the functor $\mathrm{Hom}_{\mathbf{Set}}(-,A)$.
\end{remark}
In the following paragraph we present an example of precubical set associated to a trunk.  
\paragraph{Nerve of a trunk.} We have seen in Example \ref{Example : categories and cubes are trunks} that $\square_n$ has a trunk structure and that $\partial_{i,\epsilon}^n$ is a trunk map. 
Let $\mathcal{T}$ be a trunk, the \textit{nerve of $\mathcal{T}$} is the precubical set $N(\mathcal{T}) : \square^{op} \to \mathbf{Set}$ defined by 
\begin{align}\label{precubical nerve of a trunk}
\left. 
\begin{array}{l}
N(\mathcal{T})(\square_n) := \mathrm{Hom}_{\mathbf{Trunk}}(\square_n,\mathcal{T}),\\ 
N(\mathcal{T})(\partial_{i,\epsilon}^n) := \mathrm{Hom}_{\mathbf{Trunk}}(-,\mathcal{T})(\partial_{i,\epsilon}^n) = (\partial_{i,\epsilon}^n)^*. 
\end{array}
\right.
\end{align}
In the case of the trunk associated to a rack, the nerve is easy to compute. 
\begin{theorem}\label{Theorem : nerve of a rack}
Let $X$ be a rack. There is a bijection $N(X)(\square_n) \stackrel{\eta_n}{\simeq} X^n$, and under this bijection we have
$$ 
N(X)(\partial_{i,\epsilon}^n)(x_1,\dots,x_n) = \left\{ \begin{array}{ll} 
(x_1,\dots,x_{i-1},x_i \rhd x_{i+1},\dots,x_i \rhd x_n) & \textrm{if } \epsilon = 0,\\
(x_1, \dots, x_{i-1},x_{i+1},\dots,x_n) & \textrm{if } \epsilon = 1.\\
\end{array}
\right. 
$$
\end{theorem}
\begin{lemma}\label{Lemma : nerve of a rack}
If $X$ is a rack, $F \in \mathrm{Hom}_{\mathbf{Trunk}}(\square_n,X)$ and $A \to A \amalg \{b\}$ is an edge in $\square_n$, then
\begin{align} \label{nerve of a rack}
F(A \to A \amalg \{b\}) = \displaystyle \prod_{\stackrel{1 \leq x \leq b}{x \notin A}} F([x-1] \to [x]).
\end{align}
where $[x] := \{1,\dots,x\}$.
\end{lemma}
\begin{proof}{(of Lemma \ref{Lemma : nerve of a rack})} First recall that in the trunk associated to a rack, we have $a \lhd b = a$. Suppose $A = \{a_1,\dots,a_p\}$ with $a_1 < \dots < a_p$. We can suppose that $b > a_p$. Indeed, if this is not the case, then 
\begin{align*}
F(A \to A \amalg \{b\}) &= F((A \setminus \{a_p\}) \to (A \setminus \{a_p\}) \amalg \{b\}) \lhd F((A \setminus \{a_p\}) \amalg \{b\} \to A \amalg \{b\}) \\
&= F((A \setminus \{a_p\}) \to (A \setminus \{a_p\}) \amalg \{b\})
\end{align*}
because 
$$
\xymatrix{A \ar[r] & A \amalg \{b\} \\ A \setminus \{a_p\} \ar[u] \ar[r] & (A \setminus \{a_p\}) \amalg \{b\} \ar[u]}
$$
is a preferred square in $\square_n$. Continuing this reduction until there is no element in $A$ bigger than $b$ shows that we may suppose that $b > a_p$.
\par
Now let us show by induction on the cardinal of $\{1,\dots,n\} \setminus A$ the expected equality \eqref{nerve of a rack}.
\begin{itemize}
\item{\textbf{Initialization}} Suppose that the cardinal of $\{1,\dots,n\} \setminus A$ is equal to $1$. Because $b > a_p$, necessarily $A = \{1,\dots,n-1\}$ and $b = n$, so the equality \eqref{nerve of a rack} is true.
\item{\textbf{Induction:}} Suppose the equality true for rank $k$. Let $b_1 = \min\{1 \leq x \leq b \, | \, x \notin A\}$. If $b_1 = b$, then $A = \{1,\dots,b-1\}$ and we have the expected equality. If $b_1 < b$ then the following equality holds by induction hypothesis
$$
F(A \amalg \{b_1\} \to A \amalg \{b_1,b\}) = \displaystyle \prod_{\stackrel{1 \leq x \leq b}{x \notin A \amalg \{b_1\}}} F([x-1] \to [x]).
$$ 
Because $b_1 < b$, the square 
$$
\xymatrix{A \amalg \{b\} \ar[r] & A \amalg \{b_1,b\} \\ A  \ar[u] \ar[r] & A  \amalg \{b_1\} \ar[u]}
$$
is preferred in $\square_n$, therefore
$$
F(A \to A \amalg \{b\}) = F(A \to A \amalg \{b_1\}) \rhd F(A \amalg \{b_1\} \to A \amalg \{b_1,b\}).
$$
Now using the same reduction as in the begining, we have
$$
F(A \to A \amalg \{b_1\}) = F(A_{< b_1} \to A_{< b_1} \amalg \{b_1\}),
$$
and by definition of $b_1$, we set $A_{< b_1} = \{1,\dots,b_1-1\}$. Thus the expected equality is true.
\end{itemize}
\end{proof}
\begin{proof}{(Theorem \ref{Theorem : nerve of a rack})} 
First, let $\eta_n$ the set map defined from $N(X)(\square_n)$ to $X^n$ by
\begin{align}
\eta_n : F \mapsto (x_1,\dots,x_n) \label{bijection between N(X)(squaren) and Xn}
\end{align}
with $x_{k+1} = F([k] \to [k+1])$ for all $0 \leq k \leq n-1$. By Lemma \ref{Lemma : nerve of a rack} this map is a bijection.
\par
Let $(x_1,\dots,x_n) \in X^n$. Let $F$ be the trunk map in $\mathrm{Hom}_{\mathbf{Trunk}}(\square_n,X)$ corresponding to this $n$-tuple by the previous bijection. By definition, $N(X)(\partial_{i,\epsilon}^n)(x_1,\dots,x_n) = (y_1,\dots,y_{n-1})$ with 
$$
y_k = (F \circ \partial_{i,\epsilon}^n)([k-1] \to [k]).
$$
\begin{itemize}
\item{$\epsilon = 1:$} We have
$$
\partial_{i,1}^n([k-1] \to [k])= \left\{ \begin{array}{ll} 
 [k-1] \amalg \{i\} \to  [k] \amalg \{i\} & \textrm{if } k < i,\\
 \lbrack k \rbrack \to [k+1] & \textrm{if } k \geq i.  
\end{array}
\right.
$$
Then
$$
y_k = \left\{ \begin{array}{ll} 
F([k-1] \amalg \{i\} \to [k] \amalg \{i\}) = F([k-1] \to [k]) = x_k & \textrm{if } k < i,\\
F([k] \to [k+1]) = x_{k+1} & \textrm{if } k \geq i.  
\end{array}
\right.
$$
\item{$\epsilon = 0:$} We have
$$
\partial_{i,0}^n([k-1] \to [k]) = \left\{ \begin{array}{ll} 
[k-1] \to [k]  & \textrm{if } k < i,\\
\lbrack k] \setminus \{i\} \to [k+1] \setminus \{i\} & \textrm{if } k \geq i.  
\end{array}
\right.
$$
The square
$$
\xymatrix{ [k+1] \setminus \{i\} \ar[r] & [k+1] \\ [k] \setminus \{i\} \ar[r] \ar[u] & [k] \ar[u]}
$$
is preferred in $\square_n$, so 
$$
F([k] \setminus \{i\} \to [k+1] \setminus \{i\}) = F([k] \setminus \{i\} \to [k]) \rhd F([k] \to [k+1]).
$$
By Lemma \ref{Lemma : nerve of a rack}, this is equal to 
$$
F([i-1] \to [i]) \rhd F([k] \to [k+1]),
$$
and finally
$$
y_k = \left\{ \begin{array}{ll} 
F([k-1] \to [k]) = x_k  & \textrm{if } k < i,\\
x_i \rhd x_{k+1} & \textrm{if } k \geq i.  
\end{array}
\right.
$$
\end{itemize}
\end{proof}
Let $X$ be a rack and $A$ be an abelian group. By Remark \ref{Remark : from precubical set to precocubical module}, the composition $\mathrm{Hom}_{\mathbf{Set}}(-,A) \circ N(X)$ defines a precocubical abelian group. The cochain complex associated to this precocubical abelian group is denoted by $\{C^{n}(N(X),A),d^n\}_{n \in \mathbb{N}}$.
\begin{corollary}\label{Corollary : isomorphism between CR(X,A) and C(N(X),A)}
Let $X$ be a rack and $A$ be an abelian group (considered as a trivial $X$-module). There is an isomorphism of cochain complexes
$$
\{CR^{n}(X,A),d_R^n\}_{n \in \mathbb{N}} \stackrel{\eta^*}{\simeq} \{C^{n}(N(X),A),d^n\}_{n \in \mathbb{N}}.
$$
\end{corollary}
\section{A graded dendriform algebra structure on $HR^{\bullet}(X,A)$}
Let $X$ be a rack and $A$ be an associative algebra considered as a trivial $X$-module. The goal of this section is to define a graded dendriform algebra structure on $HR^{\bullet}(X,A)$, the rack cohomology of $X$ with coefficients in $A$. First we define two graded dendriform products on $\{C^{n}(N(X),A),d^n\}_{n \in \mathbb{N}}$, the graded module associated to the nerve of the rack $X$. Then we show that these products are compatible with the differential, and so define a graded dendriform algebra structure on $HR^{\bullet}(X,A)$.
\paragraph{Products.} Let $\sigma \in \mathrm{Sh}_{p_1,p_2}$. In Example \ref{Example : functors, faces and degeneracies are trunk maps} we deduced from $\sigma$ two trunk maps $\sigma \circ i_{p_1} : \square_{p_1} \to \square_{p_1+p_2}$ and $\sigma \circ i_{p_2} : \square_{p_2} \to \square_{p_1+p_2}$. Let $\rho_{\sigma}$ denote the (set theoretical) map from $N(X)(\square_{p_1+p_2})$ to  $N(X)(\square_{p_1}) \times N(X)(\square_{p_2})$ defined by
\begin{align}
\rho_{\sigma}: = ((\sigma \circ i_{p_1})^*,(\sigma \circ i_{p_2})^*). \label{definition of rhosigma}
\end{align}
We define two graded products $\succ$ and $\prec$ on $\{C^{n}(N(X),A),d^n\}_{n \in \mathbb{N}}$ by 
\begin{align} \label{definition vdash}
f_1 \succ f_2 := \displaystyle \somme{\sigma \in \mathrm{Sh}_{p_1,p_2}^{p_1+p_2}}{} \epsilon(\sigma) \, \mu_A \circ (f_1 \times f_2) \circ \rho_{\sigma},
\end{align}
and
\begin{align} \label{definition dashv}
f_1 \prec f_2 := \displaystyle \somme{\sigma \in \mathrm{Sh}_{p_1,p_2}^{p_1}}{} \epsilon(\sigma) \, \mu_A \circ (f_1 \times f_2) \circ \rho_{\sigma}. 
\end{align}
\paragraph{Dendriform structure.} In this paragraph we prove that $\{C^n(N(X),A),\succ,\prec\}_{n \in \mathbb{N}}$ is a graded dendriform algebra. For this we use the bijections
$$
\mathrm{Sh}_{p_1,p_2+p_3} \times \mathrm{Sh}_{p_2,p_3} \stackrel{\alpha}{\simeq} \mathrm{Sh}_{p_1,p_2,p_3} \stackrel{\beta}{\simeq} \mathrm{Sh}_{p_1+p_2,p_3} \times \mathrm{Sh}_{p_1,p_2}.
$$
described in \eqref{bijections between shuffles}.
\par
First, similarly to what we did in Example \ref{Example : functors, faces and degeneracies are trunk maps}, let us define from a shuffle $\sigma \in \mathrm{Sh}_{p_1,p_2,p_3}$ three trunk maps $\sigma \circ i_{p_j} : \square_{p_j} \to \square_{p_1+p_2+p_3}$ where 
\begin{align*}
i_{p_1} &:= \partial_{p_1+p_2+p_3,0}^{p_1+p_2+p_3} \circ \dots \circ \partial_{p_1+1,0}^{p_1+1},\\
i_{p_2} &:= \partial_{p_1+p_2+p_3,0}^{p_1+p_2+p_3} \circ \dots \circ \partial_{p_1+p_2+1,0}^{p_1+p_2+1} \circ \partial_{1,1}^{p_1+p_2} \circ \dots \circ \partial_{1,1}^{p_2+1},\\
i_{p_3} &:= \partial_{1,1}^{p_1+p_2+p_3} \circ \dots \circ \partial_{1,1}^{p_3+1}.
\end{align*}
If $A = \{a_1,\dots,a_p\} \subseteq \square_{p_j}$, then 
\begin{align*}
(\sigma \circ i_{p_1})(A) &= \{\sigma(a_1),\dots,\sigma(a_p)\},\\
(\sigma \circ i_{p_2})(A) &= \{\sigma(1),\dots,\sigma(p_1),\sigma(p_1+a_1),\dots,\sigma(p_1+a_p)\},\\
(\sigma \circ i_{p_3})(A) &= \{\sigma(1),\dots,\sigma(p_1+p_2),\sigma(p_1+p_2+a_1),\dots,\sigma(p_1+p_2+a_p)\}.
\end{align*}
Let $\rho_{\sigma}$ the set map from $N(X)(\square_{p_1+p_2+p_3})$ to $N(X)(\square_{p_1}) \times N(X)(\square_{p_2}) \times N(X)(\square_{p_3})$ defined by
$$
\rho_{\sigma} :=  ((\sigma \circ i_{p_1})^*,(\sigma \circ i_{p_2})^*,(\sigma \circ i_{p_3})^*). 
$$
\begin{lemma} \label{Lemma : graded dendriform algebra structure}
Let $\sigma \in \mathrm{Sh}_{p_1,p_2+p_3}$ $(\text{resp. } \sigma \in \mathrm{Sh}_{p_1+p_2,p_3})$ and $\gamma \in \mathrm{Sh}_{p_2,p_3}$ $(\text{resp. } \gamma \in \mathrm{Sh}_{p_1,p_2})$. The following equality holds 
$$
\rho_{\sigma \circ (1 \star \gamma)} = (1 \times \rho_{\gamma}) \circ \rho_{\sigma} \, (\text{resp. } \rho_{\sigma \circ (\gamma \star 1)} = (\rho_{\gamma} \times 1) \circ \rho_{\sigma}).
$$
\end{lemma}
\begin{proof} Let $\sigma \in \mathrm{Sh}_{p_1,p_2+p_3}$ and $\gamma \in \mathrm{Sh}_{p_2,p_3}$. Let $F \in N(\mathcal{T})(\square_{p_1+p_2+p_3})$, by definition 
$$
\rho_{\sigma \circ (1 \star \gamma)}(F) = (F \circ (\sigma \circ (1 \star \gamma) \circ i_{p_1}),F \circ (\sigma \circ (1 \star \gamma) \circ i_{p_2}),F \circ (\sigma \circ (1 \star \gamma)\circ i_{p_3}))
$$ 
and 
$$
((1 \times \rho_{\gamma}) \circ \rho_{\sigma})(F) = (F \circ \sigma \circ i_{p_1},F \circ \sigma \circ i_{p_2+p_3} \circ \gamma \circ i_{p_2},F \circ \sigma_{i_{p_2+p_3}} \circ \gamma_{i_{p_3}}).
$$
Let $A = \{a_1,\dots,a_p\} \subseteq \square_{p_1}$. 
\begin{align*}
(\sigma \circ (1 \star \gamma) \circ i_{p_1})(A) &= \{(\sigma \circ (1 \star \gamma))(a_1),\dots,(\sigma \circ (1 \star \gamma))(a_p)\},\\
&=\{\sigma(a_1),\dots,\sigma(a_p)\},\\
&= (\sigma \circ i_{p_1})(A).
\end{align*}
Let $A = \{a_1,\dots,a_p\} \subseteq \square_{p_2}$.
\begin{align*}
(\sigma \circ (1 \star \gamma) \circ i_{p_2})(A) &= \{(\sigma \circ (1 \star \gamma))(1),\dots,(\sigma \circ (1 \star \gamma))(p_1),(\sigma \circ (1 \star \gamma))(p_1+a_1),\dots,(\sigma \circ (1 \star \gamma))(p_1+a_p)\},\\
&=\{\sigma(1),\dots,\sigma(p_1),\sigma(p_1+\gamma(a_1)),\dots,\sigma(p_1+\gamma(a_p))\},\\
&= (\sigma \circ i_{p_2+p_3} \circ \gamma \circ i_{p_2})(A).
\end{align*}
Let $A = \{a_1,\dots,a_p\} \subseteq \square_{p_3}$.
\begin{align*}
(\sigma \circ (1 \star \gamma)\circ i_{p_3})(A) &= \{(\sigma \circ (1 \star \gamma))(1),\dots,(\sigma \circ (1 \star \gamma))(p_1+p_2),(\sigma \circ (1 \star \gamma))(p_1+p_2+a_1),\dots\\
& \quad \dots ,(\sigma \circ (1 \star \gamma))(p_1+p_2+a_p)\},\\
&=\{\sigma(1),\dots,\sigma(p_1),\sigma(p_1+\gamma(1)),\dots,\sigma(p_1+\gamma(p_2)),\sigma(p_1+\gamma(p_2+a_1)),\dots,\sigma(p_1+\gamma(p_2+a_p))\},\\
&= (\sigma \circ i_{p_2+p_3} \circ \gamma \circ i_{p_3})(A).
\end{align*}
Thus $\rho_{\sigma \circ (1 \star \gamma)} = (1 \times \rho_{\gamma}) \circ \rho_{\sigma}$. The proof for the other equality is similar.
\end{proof}
\begin{theorem}  \label{Theorem : graded dendriform algebra structure}
Let $X$ be a rack and $A$ be an associative algebra. The graded module $\{C^{n}(N(X),A),\succ,\prec\}_{n \in \mathbb{N}}$ provided with the products $\succ$ and $\prec$ is a graded dendriform algebra.
\end{theorem}
\begin{proof}
Let $p_1,p_2,p_3 \in \mathbb{N}$ and $f_1 \in C^{p_1}(N(X),A),f_2 \in C^{p_2}(N(X),A),f_3 \in C^{p_3}(N(X),A)$. We have to prove the three equalities
\begin{enumerate}
\item $f_1 \succ (f_2 \succ f_3) = (f_1 \succ f_2 + f_1 \prec f_2) \succ f_3$,
\item $f_1 \succ (f_2 \prec f_3) = (f_1 \succ f_2) \prec f_3$,
\item $(f_1 \prec f_2) \prec f_3 = f_1 \prec (f_2 \prec f_3 + f_2 \succ f_3)$.
\end{enumerate}
Consider the bijection from $\mathrm{Sh}_{p_1,p_2+p_3} \times \mathrm{Sh}_{p_2,p_3}$ to $\mathrm{Sh}_{p_1+p_2,p_3} \times \mathrm{Sh}_{p_1,p_2}$  
$$
(\sigma,\gamma) \mapsto (\sigma',\gamma') = (\beta^{-1} \circ \alpha)(\sigma,\gamma).
$$
In particular $\epsilon(\sigma)\epsilon(\gamma) = \epsilon(\sigma')\epsilon(\gamma')$.
\begin{enumerate}
\item Let $(\sigma,\gamma) \in \mathrm{Sh}_{p_1,p_2+p_3}^{p_1+p_2+p_3} \times \mathrm{Sh}_{p_2,p_3}^{p_2+p_3}$, necessarily $(\sigma',\gamma') \in \mathrm{Sh}_{p_1+p_2,p_3}^{p_1+p_2+p_3} \times \mathrm{Sh}_{p_1,p_2}$. Then, because of cardinality, the map $\beta^{-1} \circ \alpha$ induces a bijection
$$
\mathrm{Sh}_{p_1,p_2+p_3}^{p_1+p_2+p_3} \times \mathrm{Sh}_{p_2,p_3}^{p_2+p_3} \simeq \mathrm{Sh}_{p_1+p_2,p_3}^{p_1+p_2+p_3} \times \mathrm{Sh}_{p_1,p_2}.
$$
Applying the change of variable $\beta^{-1} \circ \alpha$ to the left term in the equality $\mathrm{1.}$ and using the Lemma \ref{Lemma : graded dendriform algebra structure}, we find the right term in $\mathrm{1.}$. 
\item To show this equality, use the bijection $\beta^{-1} \circ \alpha$ restricted to $\mathrm{Sh}_{p_1,p_2+p_3}^{p_1+p_2+p_3} \times \mathrm{Sh}_{p_2,p_3}^{p_2}$. It gives a change of variables
$$
\mathrm{Sh}_{p_1,p_2+p_3}^{p_1+p_2+p_3} \times \mathrm{Sh}_{p_2,p_3}^{p_2} \simeq \mathrm{Sh}_{p_1+p_2,p_3}^{p_1+p_2} \times \mathrm{Sh}_{p_1,p_2}^{p_1+p_2}.
$$
\item To show this equality, use the bijection $\beta^{-1} \circ \alpha$ restricted to $\mathrm{Sh}_{p_1,p_2+p_3}^{p_1} \times \mathrm{Sh}_{p_1,p_2}$. It gives a change of variables
$$
\mathrm{Sh}_{p_1,p_2+p_3}^{p_1} \times \mathrm{Sh}_{p_2,p_3} \simeq \mathrm{Sh}_{p_1+p_2,p_3}^{p_1+p_2} \times \mathrm{Sh}_{p_1,p_2}^{p_1}.
$$
\end{enumerate} 
\end{proof}
\paragraph{Compatibility of the differential with the products.}
Let $X$ be a rack and $A$ be an associative algebra considered as a trivial $X$-module. By Theorem \ref{Theorem : graded dendriform algebra structure}, there is a graded dendriform algebra structure on $C^{\bullet}(N(X),A)$. In this paragraph we prove that this dendriform structure is compatible with the differential, that is, $\{C^n(N(X),A),d^n,\succ,\prec\}_{n \in \mathbb{N}}$ is a differential graded dendriform algebra.
\par
\begin{theorem}\label{Theroem : compatibility with the differential}
Let $X$ be a rack and $A$ be an associative algebra. The cochain complex $\{C^n(N(X),A),d^n\}_{n \in \mathbb{N}}$ provided with the products $\succ$ and $\prec$ is a differential graded dendriform algebra.
\end{theorem} 
\begin{proof} We have to prove that 
\begin{align}
(d^{p_1+1}f_1) \succ f_2 + (-1)^{p_1} f_1 \succ (d^{p_2+1}f_2) &= d^{p_1+p_2+1} (f_1 \succ f_2) \label{Compatibility with succ}
\end{align}
and
\begin{align}
(d^{p_1+1}f_1) \prec f_2 + (-1)^{p_1} f_1 \prec (d^{p_2+1}f_2) &= d^{p_1+p_2+1} (f_1 \prec f_2) \label{Compatibility with prec}
\end{align}
for all $f_1 \in C^{p_1}(N(X),A)$ and $f_2 \in C^{p_2}(N(X),A)$.
\vskip 0.2 cm
To prove these two equalities, let us introduce the map $\phi$ from $\mathbb{S}_{n+1} \times \{1,\dots,n+1\}$ to $\mathbb{S}_{n}$ defined by
$$
\phi(\sigma,i) := (n+1 \dots \sigma(i)) \circ \sigma \circ (i\dots n+1).
$$
Two important properties of $\phi$ are the following equalities, satisfied for all $\sigma \in \mathbb{S}_{n+1}, i \in \{1,\dots,p_1+p_2+1\}$ and $\epsilon \in \{0,1\}$ :
\begin{align}
\partial_{\sigma(i),\epsilon}^{n+1} \circ \phi(\sigma,i) &= \sigma \circ \partial_{i,\epsilon}^{n+1}, \label{change of variables 1} \\
(-1)^{\sigma(i)}\epsilon\big(\phi(\sigma,i)\big) &= (-1)^i\epsilon(\sigma). \label{change of variables 2} 
\end{align}
\begin{remark}
The equality \eqref{change of variables 1} takes place in the category $\mathbf{Gph}$.
\end{remark}
In the sequel, the restriction to suitable subsets of the map $\psi :\mathbb{S}_{n+1} \times \{1,\dots,n+1\} \to \mathbb{S}_{n} \times \{1,\dots,n+1\}$ defined by $\psi(\sigma,i) = (\phi(\sigma,i),\sigma(i))$ will provide the right change of variables needed in order to prove \eqref{Compatibility with succ} and \eqref{Compatibility with prec}.
\vskip 0.2 cm
\noindent
\textit{Compatibility with $\succ$}: Let $f_1 \in C^{p_1}(N(X),A)$ and $f_2 \in C^{p_2}(N(X),A)$. The left term in \eqref{Compatibility with succ} is equal to the sum of four terms $\mathrm{(a), (b), (c)}$ and $\mathrm{(d)}$ with
\begin{align*}
\displaystyle
\mathrm{(a)} &= \sum_{\sigma \in \mathrm{Sh}_{p_1+1,p_2}^{p_1+p_2+1}}^{} \sum_{i = 1}^{p_1} (-1)^i \epsilon(\sigma) \, \mu_A \circ (f_1 \times f_2) \circ \big((\sigma \circ i_{p_1+1} \circ \partial_{i,0}^{p_1+1})^*,(\sigma \circ i_{p_2})^*\big), \\ 
\mathrm{(b)} &= \sum_{\sigma \in \mathrm{Sh}_{p_1,p_2+1}^{p_1+p_2+1}} \sum_{i = 1}^{p_2} (-1)^{i + p_1} \epsilon(\sigma) \, \mu_A \circ (f_1 \times f_2) \circ \big((\sigma \circ i_{p_1})^*,(\sigma \circ i_{p_2+1}  \circ \partial_{i,0}^{p_2+1})^*\big), \\
\mathrm{(c)} &= \sum_{\sigma \in \mathrm{Sh}_{p_1+1,p_2}^{p_1+p_2+1}}^{} \sum_{i = 1}^{p_1} (-1)^{i+1} \epsilon(\sigma) \, \mu_A \circ (f_1 \times f_2) \circ \big((\sigma \circ i_{p_1+1} \circ \partial_{i,1}^{p_1+1})^*,(\sigma \circ i_{p_2})^*\big), \\ 
\mathrm{(d)} &= \sum_{\sigma \in \mathrm{Sh}_{p_1,p_2+1}^{p_1+p_2+1}} \sum_{i = 1}^{p_2} (-1)^{i + p_1+1} \epsilon(\sigma) \, \mu_A \circ (f_1 \times f_2) \circ \big((\sigma \circ i_{p_1})^*,(\sigma \circ i_{p_2+1}  \circ \partial_{i,1}^{p_2+1})^*\big). \\  
\end{align*}
The right term in \eqref{Compatibility with succ} is equal to the sum of two terms $\mathrm{(A)}$ and $\mathrm{(B)}$, with
\begin{align*}
\displaystyle
\mathrm{(A)} &= \sum_{\sigma \in \mathrm{Sh}_{p_1,p_2}^{p_1+p_2}}^{} \sum_{i = 1}^{p_1+p_2} (-1)^{i} \epsilon(\sigma) \, \mu_A \circ (f_1 \times f_2) \circ \big((\partial_{i,0}^{p_1+p_2+1} \circ \sigma \circ i_{p_1})^*, (\partial_{i,0}^{p_1+p_2+1} \circ \sigma \circ i_{p_2})^*\big),\\
\mathrm{(B)} &= \sum_{\sigma \in \mathrm{Sh}_{p_1,p_2}^{p_1+p_2}}^{} \sum_{i = 1}^{p_1+p_2} (-1)^{i+1} \epsilon(\sigma) \, \mu_A \circ (f_1 \times f_2) \circ \big((\partial_{i,1}^{p_1+p_2+1} \circ \sigma \circ i_{p_1})^*, (\partial_{i,1}^{p_1+p_2+1} \circ \sigma \circ i_{p_2})^*\big).
\end{align*}
In order to prove \eqref{Compatibility with succ}, we are going to show that
$$
\mathrm{(a)}+\mathrm{(c)}+\mathrm{(d)} = \mathrm{(B)} \, \text{ and } \, \mathrm{(b)} = \mathrm{(A)}.
$$
\vskip 0.3 cm
\noindent
$\bullet \mathrm{(b)} = \mathrm{(A)}:$ This equality is proved using the change of variables
$$
\psi \circ t_{p_1} : \mathrm{Sh}_{p_1,p_2+1}^{p_1+p_2+1} \times \{1,\dots,p_2\} \to \mathrm{Sh}_{p_1,p_2}^{p_1+p_2} \times \{1,\dots,p_1+p_2\}
$$ 
with $t_{p_1}(i) = p_1+i$. Using this change of variables, we have to prove for all $\sigma \in \mathrm{Sh}_{p_1,p_2+1}^{p_1+p_2+1}$ and $i \in \{1,\dots,p_2\}$ :
\begin{enumerate}
\item $\big((\sigma \circ i_{p_1})^*,(\sigma \circ i_{p_2+1}  \circ \partial_{i,0}^{p_2+1})^*\big) = \big((\partial_{\sigma(p_1+i),0}^{p_1+p_2+1} \circ \phi(\sigma,p_1+i) \circ i_{p_1})^*, (\partial_{\sigma(p_1+i),0}^{p_1+p_2+1} \circ \phi(\sigma,p_1+i) \circ i_{p_2})^*\big)$,
\item $(-1)^{\sigma(p_1+i)} \epsilon\big(\phi(\sigma,p_1+i)\big) = (-1)^{p_1+i}\epsilon(\sigma)$. 
\end{enumerate}
\vskip 0.3 cm
\noindent
$\mathrm{1.}$ Let $\sigma \in \mathrm{Sh}_{p_1,p_2+1}^{p_1+p_2+1}$ and $i \in \{1,\dots,p_2\}$. We want to prove for all $F \in \mathrm{Hom}_{\mathbf{Trunk}}(\square_{p_1+p_2+1},X)$
\begin{align*}
\left\{\begin{array}{l}
F \circ \sigma \circ i_{p_1} =  F \circ \partial_{\sigma(p_1+i),0}^{p_1+p_2+1} \circ \phi(\sigma,p_1+i) \circ i_{p_1}, \\ 
F \circ \sigma \circ i_{p_2+1} \circ \partial_{i,0}^{p_2+1} = F \circ \partial_{\sigma(p_1+i),0}^{p_1+p_2+1} \circ \phi(\sigma,p_1+i) \circ i_{p_2}.  
\end{array}\right.
\end{align*}
By Theorem \ref{Theorem : nerve of a rack}, this is equivalent to
\begin{align*}
\left\{\begin{array}{ll}
(F \circ \sigma \circ i_{p_1})([k-1] \to [k]) =  (F \circ \partial_{\sigma(p_1+i),0}^{p_1+p_2+1} \circ \phi(\sigma,p_1+i) \circ i_{p_1})([k-1] \to [k]) & 1 \leq k \leq p_1, \\ 
(F \circ \sigma \circ i_{p_2+1} \circ \partial_{i,0}^{p_2+1})([k-1] \to [k]) = (F \circ \partial_{\sigma(p_1+i),0}^{p_1+p_2+1} \circ \phi(\sigma,p_1+i) \circ i_{p_2})([k-1] \to [k]) & 1 \leq k \leq p_2.  
\end{array}\right.
\end{align*}
Let $0 \leq k \leq p_1$, $(\sigma \circ i_{p_1})([k]) = \sigma([k])$ and $(\partial_{\sigma(p_1+i),0}^{p_1+p_2+1} \circ \phi(\sigma,p_1+i) \circ i_{p_1})([k]) \stackrel{\eqref{change of variables 1}}{=} (\sigma \circ \partial_{p_1+i,0}^{p_1+p_2+1} \circ i_{p_1})([k]) = \sigma([k])$, so the first equality is true.\\
Let $0 \leq k \leq p_2$,  
\begin{align*}
(\sigma \circ i_{p_2+1} \circ \partial_{i,0}^{p_2+1})([k]) &= \left\{\begin{array}{ll}
										(\sigma \circ i_{p_2+1})([k]) & \textrm{if } k < i,\\
										(\sigma \circ i_{p_2+1})([k+1] \setminus \{i\}) & \textrm{if } k \geq i,
										\end{array}\right.\\
										&= \left\{\begin{array}{ll}
										\sigma([p_1+k]) & \textrm{if } k < i,\\
										\sigma([p_1+k+1] \setminus \{p_1+i\}) & \textrm{if } k \geq i,
										\end{array}\right.
\end{align*}
and 
\begin{align*}
\partial_{\sigma(p_1+i),0}^{p_1+p_2+1} \circ \phi(\sigma,p_1+i) \circ i_{p_2})([k]) &= \partial_{\sigma(p_1+i),0}^{p_1+p_2+1} \circ \phi(\sigma,p_1+i))([p_1+k]) \\
& \stackrel{\eqref{change of variables 1}}{=} (\sigma \circ \partial_{p_1+i,0}^{p_1+p_2+1})([p_1+k]) \\
& =  \left\{\begin{array}{ll}
										\sigma([p_1+k]) & \textrm{if } k < i,\\
										\sigma([p_1+k+1] \setminus \{p_1+i\}) & \textrm{if } k \geq i,
										\end{array}\right.
\end{align*}  
so the second equality is true.
\vskip 0.3 cm
\noindent
$\mathrm{2.}$ Deduced from \eqref{change of variables 2}.
\vskip 0.3 cm
\noindent
$\bullet \mathrm{(a)+(c)+(d)} = \mathrm{(B)}:$ To prove this equality, we need to decompose $\mathrm{(d)}$ into two sums. Let us denote by $I^{\geq}$ and $I^{<}$ the subsets of $\mathrm{Sh}_{p_1,p_2+1}^{p_1+p_2+1} \times \{1,\dots,p_2\}$ defined by
\begin{align*}
I^{\geq} &:= \{(\sigma,i) \in \mathrm{Sh}_{p_1,p_2+1}^{p_1+p_2+1} \times \{1,\dots,p_2\} \, | \, \sigma(p_1+i) \geq p_1+i\} \\
& \, \, = \{(\sigma,i) \in \mathrm{Sh}_{p_1,p_2+1}^{p_1+p_2+1} \times \{1,\dots,p_2\} \, | \, \sigma(p_1+i) = p_1+i\}  
\end{align*}
and
$$
I^{<} := (\mathrm{Sh}_{p_1,p_2+1}^{p_1+p_2+1} \times \{1,\dots,p_2\}) \setminus I^{\geq}.
$$
Then we are going to prove that 
$$
\displaystyle 
\mathrm{(a)} + \sum_{I^{<}}^{} = 0 \, \text{ and } \, \displaystyle \mathrm{(c)} + \sum_{I^{\geq}}^{} = \mathrm{(B)}.
$$
$\mathrm{(a)} + \sum_{I^{<}}^{} = 0:$ To prove this equality, we use the change of variables
$$
\theta : \mathrm{Sh}_{p_1+1,p_2}^{p_1+p_2+1} \times \{1,\dots,p_1\} \to I^{<}
$$
defined by $\theta(\sigma,i) = (\rho(\sigma,i),1+\sigma(i)-i)$ with $\rho(\sigma,i) = \sigma \circ (i \dots p_1+1+\sigma(i)-i)$. Using this change of variables, we have to prove for all $\sigma \in \mathrm{Sh}_{p_1+1,p_2}^{p_1+p_2+1}$ and $i \in \{1,\dots,p_1\}$ :
\begin{enumerate}
\item $\big((\sigma \circ i_{p_1+1} \circ \partial_{i,0}^{p_1+1})^*,(\sigma \circ i_{p_2})^*\big) = \big((\rho(\sigma,i) \circ i_{p_1})^*,(\rho(\sigma,i) \circ i_{p_2+1}  \circ \partial_{1+\sigma(i)-i,1}^{p_2+1})^*\big)$,
\item $(-1)^{p_1+1+\sigma(i)-i}\epsilon\big(\rho(\sigma,i)\big) = (-1)^{i}\epsilon(\sigma)$.
\end{enumerate}
\vskip 0.3 cm
\noindent
$\mathrm{1.}$ Let $\sigma \in \mathrm{Sh}_{p_1+1,p_2}^{p_1+p_2+1}$ and $i \in \{1,\dots,p_1\}$. Using Theorem \ref{Theorem : nerve of a rack}, we are going to prove, for all $F \in \mathrm{Hom}_{\mathbf{Trunk}}(\square_{p_1+p_2+1},X)$, the following equalities :
$$
\left\{ \begin{array}{ll}
(F \circ \sigma \circ i_{p_1+1} \circ \partial_{i,0}^{p_1+1})([k-1] \to [k]) = (F \circ \rho(\sigma,i) \circ i_{p_1})([k-1] \to [k]) & \textrm{if } 1 \leq k \leq p_1,\\
(F \circ \sigma \circ i_{p_2})([k-1] \to [k]) = (F \circ \rho(\sigma,i) \circ i_{p_2+1} \circ \partial_{1 + \sigma(i)-i,1}^{p_2+1})([k-1] \to [k]) & \textrm{if } 1 \leq k \leq p_2.
\end{array}\right.
$$
Let $0 \leq k \leq p_1$,
\begin{align*}
(\sigma \circ i_{p_1+1} \circ \partial_{i,0}^{p_1+1})([k]) &= \left\{\begin{array}{ll}
										\sigma([k]) & \textrm{if } k < i,\\
										\sigma([k+1] \setminus \{i\}) & \textrm{if } k \geq i.
										\end{array}\right.
\end{align*}
and
\begin{align*}
(\rho(\sigma,i) \circ i_{p_1})([k]) &= \big(\sigma \circ (i \dots p_1+1+\sigma(i)-i)\big)([k]),\\
						&= \left\{\begin{array}{ll}
										\sigma([k]) & \textrm{if } k < i,\\
										\sigma([k+1] \setminus \{i\}) & \textrm{if } k \geq i,
										\end{array}\right.
\end{align*}
(remark that $\sigma(i)-i \geq 0$ for $i \in \{1,\dots,p_1\}$), so the first equality is true.
\\
Let $0 \leq k \leq p_2$, we have $(\sigma \circ i_{p_2})([k]) = \sigma([p_1+1+k])$ and
\begin{align*}
(\rho(\sigma,i) \circ i_{p_2+1} \circ \partial_{1+\sigma(i)-i,1}^{p_2+1})([k]) &= \left\{\begin{array}{ll}
								(\rho(\sigma,i) \circ i_{p_2+1})([k] \amalg \{1+\sigma(i)-i\} & \textrm{if } k < 1+\sigma(i)-i,\\
								(\rho(\sigma,i) \circ i_{p_2+1})([k+1]) & \textrm{if } k \geq 1+\sigma(i)-i,\\
								\end{array}\right.\\
								&= \left\{\begin{array}{ll}
								\rho(\sigma,i)([p_1+k] \amalg \{p_1+1+\sigma(i)-i\}) & \textrm{if } k < 1+\sigma(i)-i,\\
								\rho(\sigma,i)([p_1+k+1]) & \textrm{if } k \geq 1+\sigma(i)-i,\\
								\end{array}\right.\\
								&= \left\{\begin{array}{ll}
								\sigma([p_1+k+1]) \amalg \{\sigma(i)\} & \textrm{if } k < 1+\sigma(i)-i,\\
								\sigma([p_1+k+1]) & \textrm{if } k \geq 1+\sigma(i)-i,\\
								\end{array}\right.\\
								&= \sigma([p_1+k+1]),
\end{align*}
so the second equality is true.
\vskip 0.3 cm
\noindent
$\mathrm{2.}$ Clear by definition of $\rho(\sigma,i)$.
\vskip 0.3cm
\noindent
$\displaystyle \mathrm{(c)} + \sum_{I^{\geq}}^{} = \mathrm{(B)}:$ To prove this equality, we need to decompose $\mathrm{(B)}$ into two sums, one equals $\sum_{I^{\geq}}^{}$ and the other equals $\mathrm{(c)}$. Let us denote by $J^{\geq}$ and $J^{<}$ the subsets of $\mathrm{Sh}_{p_1,p_2}^{p_1+p_2} \times \{p_1+1,\dots,p_1+p_2\}$ defined by
$$
J^{\geq} := \{(\sigma,i) \in \mathrm{Sh}_{p_1,p_2}^{p_1+p_2} \times \{p_1+1,\dots,p_1+p_2\} \, | \, \sigma(i) \geq i\}
$$ 
and
$$
J^{<} := (\mathrm{Sh}_{p_1,p_2}^{p_1+p_2} \times \{p_1+1,\dots,p_1+p_2\}) \setminus J^{\geq}.
$$
Then we are going to prove that 
$$
\displaystyle 
\sum_{I^{\geq}}^{} = \sum_{J^{\geq}}^{} \, \text{ and } \, \mathrm{(c)} = \sum_{J^{<}}^{}.
$$
\vskip 0.3 cm
\noindent
$\sum_{I^{\geq}}^{} = \sum_{J^{\geq}}^{}:$ This equality is proved using the change of variables
$$
\psi \circ t_{p_1} : I^{\geq} \to J^{\geq}.
$$
In this case we have $(\psi \circ t_{p_1})(\sigma,i) = (\sigma,p_1+i)$. Using this change of variables, we have to prove for all $\sigma \in \mathrm{Sh}_{p_1,p_2+1}^{p_1+p_2+1}$ and $i \in \{1,\dots,p_2\}$ :
$$
\big((\sigma \circ i_{p_1})^*,(\sigma \circ i_{p_2+1}  \circ \partial_{i,1}^{p_2+1})^*\big) = \big((\partial_{p_1+i,1}^{p_1+p_2+1} \circ \sigma \circ i_{p_1})^*, (\partial_{p_1+i,1}^{p_1+p_2+1} \circ \sigma \circ i_{p_2})^*\big).
$$
Let $\sigma \in \mathrm{Sh}_{p_1,p_2+1}^{p_1+p_2+1}$ and $i \in \{1,\dots,p_2\}$. Using Theorem \ref{Theorem : nerve of a rack}, we are going to prove, for all $F \in \mathrm{Hom}_{\mathbf{Trunk}}(\square_{p_1+p_2+1},X)$, the following equalities :
$$
\left\{
\begin{array}{ll}
(F \circ \sigma \circ i_{p_1})([k-1] \to [k]) = (F \circ \partial_{p_1+i,1}^{p_1+p_2+1} \circ \sigma \circ i_{p_1})([k-1] \to [k]) & \textrm{if } 1 \leq k \leq p_1,\\
(F \circ \sigma \circ i_{p_2+1}  \circ \partial_{i,1}^{p_2+1})([k-1] \to [k]) = (F \circ \partial_{p_1+i,1}^{p_1+p_2+1} \circ \sigma \circ i_{p_2})([k-1] \to [k]) & \textrm{if } 1 \leq k \leq p_2.
\end{array}
\right.
$$
Let $0 \leq k \leq p_1$, $(\sigma \circ i_{p_1})([k]) = \sigma([k])$ and $(\partial_{p_1+i,1}^{p_1+p_2+1} \circ \sigma \circ i_{p_1})([k]) \stackrel{\eqref{change of variables 1}}{=} (\sigma \circ \partial_{p_1+i,1}^{p_1+p_2+1} \circ i_{p_1})([k]) = \sigma([k]) \amalg \{\sigma(p_1+i)\} =  \sigma([k]) \amalg \{p_1+i\}$. Hence the first equality is equivalent to
$$
F(\sigma([k-1]) \to \sigma([k])) = F\big((\sigma([k-1]) \amalg \{p_1+i\}) \to (\sigma([k]) \amalg \{p_1+i\})\big).
$$
Because $\sigma \in I^{\geq}$, necessarily $\sigma(k) < p_1+i$, so the square
$$
\xymatrix{
\sigma([k-1]) \amalg \{p_1+i\} \ar[r] & \sigma([k]) \amalg \{p_1+i\} \\
\sigma([k-1]) \ar[u] \ar[r] & \sigma([k]) \ar[u]
}
$$
is a preferred square in $\square_{p_1+p_2+1}$, and the expected equality si proved.\\
Let $0 \leq k \leq p_2$, we have 
\begin{align*}
(\sigma \circ i_{p_2+1} \circ \partial_{i,1}^{p_2+1})([k]) &= \left\{\begin{array}{ll}
	\sigma([p_1+k]) \amalg \{p_1+i\} & \textrm{if } k < i,\\
	\sigma([p_1+k+1]) & \textrm{if } k \geq i,
	\end{array}
	\right.
\end{align*}
and
\begin{align*}
(\partial_{p_1+i,1}^{p_1+p_2+1} \circ \sigma \circ i_{p_2})([k]) &\stackrel{\eqref{change of variables 1}}{=} (\sigma \circ \partial_{p_1+i,1}^{p_1+p_2+1} \circ i_{p_2})([k]),\\
&= \left\{\begin{array}{ll}
\sigma([p_1+k]) \amalg \{p_1+i\} & \textrm{if } k < i,\\
\sigma([p_1+k+1]) & \textrm{if } k \geq i,
\end{array}
\right.
\end{align*}
so the second equality is true.
\vskip 0.3 cm
\noindent
$\mathrm{(c)} = \sum_{J^{<}}^{}:$ The following equality is proved using the change of variables
$$
\psi : \mathrm{Sh}_{p_1+1,p_2}^{p_1+p_2+1} \times \{1,\dots,p_1\} \to J^{<}.
$$
Using this change of variables, we have to prove for all $\sigma \in \mathrm{Sh}_{p_1+1,p_2}^{p_1+p_2+1}$ and $i \in \{1,\dots,p_1\}$ :
\begin{enumerate}
\item $\big((\sigma \circ i_{p_1+1} \circ \partial_{i,1}^{p_1+1})^*,(\sigma \circ i_{p_2})^*\big) = \big((\partial_{\sigma(i),1}^{p_1+p_2+1} \circ \phi(\sigma,i) \circ i_{p_1})^*, (\partial_{\sigma(i),1}^{p_1+p_2+1} \circ \phi(\sigma,i) \circ i_{p_2})^*\big)$,
\item $(-1)^{i+1}\epsilon(\sigma) = (-1)^{\sigma(i)+1}\epsilon(\phi(\sigma,i))$.
\end{enumerate}
$\mathrm{1.}$ Let $\sigma \in \mathrm{Sh}_{p_1+1,p_2}^{p_1+p_2+1}$ and $i \in \{1,\dots,p_1\}$. Using Theorem \ref{Theorem : nerve of a rack}, we are going to prove, for all $F \in \mathrm{Hom}_{\mathbf{Trunk}}(\square_{p_1+p_2+1},X)$, the following equalities :
$$
\left\{
\begin{array}{ll}
(F \circ \sigma \circ i_{p_1+1} \circ \partial_{i,1}^{p_1+1})([k-1] \to [k]) = (F \circ \partial_{\sigma(i),1}^{p_1+p_2+1} \circ \phi(\sigma,i) \circ i_{p_1})([k-1] \to [k]) & \textrm{if } 1 \leq k \leq p_1,\\
(F \circ \sigma \circ i_{p_2})([k-1] \to [k]) = (F \circ \partial_{\sigma(i),1}^{p_1+p_2+1} \circ \phi(\sigma,i) \circ i_{p_2})([k-1] \to [k]) & \textrm{if } 1 \leq k \leq p_2.
\end{array}
\right.
$$
Let $0 \leq k \leq p_1$, we have 
$$
(\sigma \circ i_{p_1+1} \circ \partial_{i,1}^{p_1+1})([k]) = \left\{\begin{array}{ll}
\sigma([k]) \amalg \{\sigma(i)\} & \textrm{if } k < i,\\
\sigma([k+1]) & \textrm{if } k \geq i,
\end{array}
\right.
$$
and
\begin{align*}
(\partial_{\sigma(i),1}^{p_1+p_2+1} \circ \phi(\sigma,i) \circ i_{p_1})([k]) &\stackrel{\eqref{change of variables 1}}{=} (\sigma \circ \partial_{i,1}^{p_1+p_2+1} \circ i_{p_1})([k]), \\
&= \left\{\begin{array}{ll}
\sigma([k]) \amalg \{\sigma(i)\} & \textrm{if } k < i,\\
\sigma([k+1]) & \textrm{if } k \geq i,
\end{array} 
\right. 
\end{align*}
so the first equality is true.\\
Let $0 \leq k \leq p_2$, we have $(\sigma \circ i_{p_2})([k]) = \sigma([p_1+1+k])$ and $(\partial_{\sigma(i),1}^{p_1+p_2+1} \circ \phi(\sigma,i) \circ i_{p_2})([k]) \stackrel{\eqref{change of variables 1}}{=} (\sigma \circ \partial_{i,1}^{p_1+p_2+1} \circ i_{p_2})([k]) = (\sigma \circ \partial_{i,1}^{p_1+p_2+1})([p_1+k]) = \sigma([p_1+1+k])$, so the second equality is true.
\vskip 0.3 cm
\noindent
$\mathrm{2.}$ Deduced from \eqref{change of variables 2}.
\vskip 0.3cm
\noindent
\textit{Compatibility with $\prec$:} The proof is essentially the same as before. Let $f_1 \in C^{p_1}(N(X),A)$ and $f_2 \in C^{p_2}(N(X),A)$. The left term in \eqref{Compatibility with prec} is equal to the sum of four terms $\mathrm{(a), (b), (c)}$ and $\mathrm{(d)}$ with
\begin{align*}
\displaystyle
\mathrm{(a)} &= \sum_{\sigma \in \mathrm{Sh}_{p_1+1,p_2}^{p_1+1}}^{} \sum_{i = 1}^{p_1} (-1)^i \epsilon(\sigma) \, \mu_A \circ (f_1 \times f_2) \circ \big((\sigma \circ i_{p_1+1} \circ \partial_{i,0}^{p_1+1})^*,(\sigma \circ i_{p_2})^*\big), \\ 
\mathrm{(b)} &= \sum_{\sigma \in \mathrm{Sh}_{p_1,p_2+1}^{p_1}} \sum_{i = 1}^{p_2} (-1)^{i + p_1} \epsilon(\sigma) \, \mu_A \circ (f_1 \times f_2) \circ \big((\sigma \circ i_{p_1})^*,(\sigma \circ i_{p_2+1}  \circ \partial_{i,0}^{p_2+1})^*\big), \\
\mathrm{(c)} &= \sum_{\sigma \in \mathrm{Sh}_{p_1+1,p_2}^{p_1+1}}^{} \sum_{i = 1}^{p_1} (-1)^{i+1} \epsilon(\sigma) \, \mu_A \circ (f_1 \times f_2) \circ \big((\sigma \circ i_{p_1+1} \circ \partial_{i,1}^{p_1+1})^*,(\sigma \circ i_{p_2})^*\big), \\ 
\mathrm{(d)} &= \sum_{\sigma \in \mathrm{Sh}_{p_1,p_2+1}^{p_1}} \sum_{i = 1}^{p_2} (-1)^{i + p_1+1} \epsilon(\sigma) \, \mu_A \circ (f_1 \times f_2) \circ \big((\sigma \circ i_{p_1})^*,(\sigma \circ i_{p_2+1}  \circ \partial_{i,1}^{p_2+1})^*\big). \\  
\end{align*}
The right term in \eqref{Compatibility with prec} is equal to the sum of two terms $\mathrm{(A)}$ and $\mathrm{(B)}$, with
\begin{align*}
\displaystyle
\mathrm{(A)} &= \sum_{\sigma \in \mathrm{Sh}_{p_1,p_2}^{p_1}}^{} \sum_{i = 1}^{p_1+p_2} (-1)^{i} \epsilon(\sigma) \, \mu_A \circ (f_1 \times f_2) \circ \big((\partial_{i,0}^{p_1+p_2+1} \circ \sigma \circ i_{p_1})^*, (\partial_{i,0}^{p_1+p_2+1} \circ \sigma i_{p_2})^*\big),\\
\mathrm{(B)} &= \sum_{\sigma \in \mathrm{Sh}_{p_1,p_2}^{p_1}}^{} \sum_{i = 1}^{p_1+p_2} (-1)^{i+1} \epsilon(\sigma) \, \mu_A \circ (f_1 \times f_2) \circ \big((\partial_{i,1}^{p_1+p_2+1} \circ \sigma \circ i_{p_1})^*, (\partial_{i,1}^{p_1+p_2+1} \circ \sigma \circ i_{p_2})^*\big).
\end{align*}
In order to prove \eqref{Compatibility with prec}, we are going to show that
$$
\mathrm{(a)}+\mathrm{(b)}+\mathrm{(d)} = \mathrm{(A)} \, \text{ and } \, \mathrm{(c)} = \mathrm{(B)}.
$$
\vskip 0.2cm
\noindent
$\bullet \mathrm{(c) = (B)}:$ The following equality is proved using the change of variables
$$
\psi : \mathrm{Sh}_{p_1+1,p_2}^{p_1+1} \times \{1,\dots,p_1\} \to \mathrm{Sh}_{p_1,p_2}^{p_1} \times \{1,\dots,p_1+p_2\}.
$$
\vskip 0.2cm
\noindent
$\bullet \mathrm{(a)+(b)+(d) = (A)}:$ To prove this equality, we need to decompose (a) into two sums. Let us denote by $I^{\geq}$ and $I^{<}$ the subsets of $\mathrm{Sh}_{p_1+1,p_2}^{p_1+1} \times \{1,\dots,p_1\}$ defined by
\begin{align*}
I^{\geq} &:= \{(\sigma,i) \in \mathrm{Sh}_{p_1+1,p_2}^{p_1+1} \times \{1,\dots,p_1\} \, | \, \sigma(i) \geq p_2+i\},\\
&= \{(\sigma,i) \in \mathrm{Sh}_{p_1+1,p_2}^{p_1+1} \times \{1,\dots,p_1\} \, | \, \sigma(i) = p_2+i\}.
\end{align*}
and
$$
I^{<} := (\mathrm{Sh}_{p_1+1,p_2}^{p_1+1} \times \{1,\dots,p_1\}) \setminus I^{<}.
$$
Then we are going to show that
$$\displaystyle
\mathrm{(d)} + \sum_{I^{<}}^{} = 0 \, \text{ and } \mathrm{(b)} + \sum_{I^{\geq}}^{} = \mathrm{(A)}.
$$
\vskip 0.2cm
\noindent
$\mathrm{(d)} + \sum_{I^{<}}^{} = 0:$ The following equality is proved using the change of variables
$$
\theta : \mathrm{Sh}_{p_1,p_2+1}^{p_1} \times \{1,\dots,p_2\} \to I^{<}
$$
defined by $\theta(\sigma,i) = (\rho(\sigma,i),\sigma(p_1+i)-i+1)$ with $\rho(\sigma,i) = \sigma \circ (p_1+i \, p_1+i-1 \dots \gamma(p_1+i)-i+1)$.
\vskip 0.2cm
\noindent
$\displaystyle \mathrm{(b)} + \sum_{I^{\geq}}^{} = \mathrm{(A)}:$ To prove this equality, we need to decompose (A) into two sums, one equals to $\sum_{I^{\geq}}^{}$ and the other equals to (b). Let us denote by $J^{\geq}$ and $J^{<}$ the subsets of $\mathrm{Sh}_{p_1,p_2}^{p_1} \times \{1,\dots,p_1\}$ defined by
$$
J^{\geq} := \{(\sigma,i) \in \mathrm{Sh}_{p_1,p_2}^{p_1} \times \{p_2+1,\dots,p_1+p_2\} \, | \, \sigma(i-p_2) \geq i\}
$$  
and
$$
J^{<} := (\mathrm{Sh}_{p_1,p_2}^{p_1} \times \{1,\dots,p_1\})  \setminus J^{\geq}.
$$
Then we are going to show that
$$\displaystyle
\sum_{I^{\geq}}^{} = \sum_{I^{\geq}}^{} \, \text{ and } \, \mathrm{(b)} = \sum_{J^<}^{}.
$$
\vskip0.2cm
\noindent
$\sum_{I^{\geq}}^{} = \sum_{J^{\geq}}^{}:$ The following equality is proved using the change of variables
$$
\psi : I^{\geq} \to J^{\geq}.
$$
\vskip0.2cm
\noindent
$\mathrm{(b)} = \sum_{J^{<}}^{}:$ The following equality is proved using the change of variables
$$
\psi \circ t_{p_1} : \mathrm{Sh}_{p_1,p_2+1}^{p_1} \times \{1,\dots,p_2\} \to J^{<}.
$$
\end{proof}
\paragraph{Explicit formulas for $\succ$ and $\prec$.} Using the cochain complex isomorphism $\eta^{*}$ between $\{CR^n(X,A),d_R^n\}_{n \in \mathbb{N}}$ and $\{C^n(N(X),A),d^n\}_{n \in \mathbb{N}}$ (Corollary \ref{Corollary : isomorphism between CR(X,A) and C(N(X),A)}), the differential graded dendriform algebra structure on $\{CR^n(X,A),d_R^n\}_{n \in \mathbb{N}}$ is defined by the formulas :
$$
f_1 \succ f_2 := \sum_{\sigma \in \mathrm{Sh}_{p_1,p_2}^{p_1+p_2}}^{} \epsilon(\sigma) \, \mu_A \circ (f_1 \times f_2) \circ (\eta_{p_1} \times \eta_{p_2}) \circ \rho_{\sigma} \circ \eta_{p_1+p_2}^{-1}, 
$$
and
$$
f_1 \prec f_2 := \sum_{\sigma \in \mathrm{Sh}_{p_1,p_2}^{p_1}}^{} \epsilon(\sigma) \, \mu_A \circ (f_1 \times f_2) \circ (\eta_{p_1} \times \eta_{p_2}) \circ \rho_{\sigma} \circ \eta_{p_1+p_2}^{-1}.
$$
Then to have formulas for $\succ$ and $\prec$ on $\{CR^n(X,A),d_R^n\}$ we have to compute 
$$
(\eta_{p_1} \times \eta_{p_2}) \circ \rho_{\sigma} \circ \eta_{p_1+p_2}^{-1} : X^{p_1+p_2} \to X^{p_1} \times X^{p_2}.
$$
\par
Let $(x_1,\dots,x_{p_1+p_2}) \in X^{p_1+p_2}$, by \eqref{bijection between N(X)(squaren) and Xn} and \eqref{definition of rhosigma}
$$
\big((\eta_{p_1} \times \eta_{p_2}) \circ \rho_{\sigma} \circ \eta_{p_1+p_2}^{-1}\big)(x_1,\dots,x_{p_1+p_2}) = \big((y_1,\dots,y_{p_1}),(z_{1},\dots,z_{p_2})\big),
$$
where 
$$
y_{k} = (F \circ \sigma \circ i_{p_1})([k-1] \to [k]) \, \text{ and } \, z_k = (F \circ \sigma \circ i_{p_2})([k-1] \to [k]).
$$
Thus, using Lemma \eqref{Lemma : nerve of a rack}
\begin{align*}
y_k &= (F \circ \sigma)([k-1] \to [k]),\\
&= F\big(\sigma([k-1]) \to \sigma([k])\big),\\
&= F\big(\sigma([k-1]) \to \sigma([k-1] \amalg \{\sigma(k)\})\big),\\
&= \displaystyle \prod_{\stackrel{1 \leq x \leq \sigma(k)}{x \notin \sigma([k-1])}} F([x-1] \to [x]),\\
&= x_{i_1} \rhd \dots \rhd x_{i_j} \rhd x_{\sigma(k)},
\end{align*}
with $x_{i_1} < \dots < x_{i_j} < \sigma(k)$ and $i_l \in \{\sigma(p_1+1),\dots,\sigma(p_1+p_2)\}$. In the same way
\begin{align*}
z_k &= (F \circ \sigma)([p_1+k-1] \to [p_1+k]),\\
&= F\big(\sigma([p_1+k-1]) \to \sigma([p_1+k])\big),\\
&= F\big(\sigma([p_1+k-1]) \to \sigma([p_1+k-1] \amalg \{\sigma(p_1+k)\})\big),\\
&= \displaystyle \prod_{\stackrel{1 \leq x \leq \sigma(p_1+k)}{x \notin \sigma([p_1+k-1])}} F([x-1] \to [x]),\\
&= x_{\sigma(p_1+k)}.
\end{align*}
Finally we have proved the following theorem.
\begin{theorem}\label{Theorem : CR(X,A) is a differential graded dendriform algebra}
Let $X$ be a rack and $A$ be an associative algebra. The cochain complex $\{CR^n(X,A),d_R^n\}_{n \in \mathbb{N}}$ provided with the products $\succ$ and $\prec$ defined by 
\begin{align*}
(f_1 \succ f_2)(x_1,\dots,x_{p_1+p_2}) &= \sum_{\sigma \in \mathrm{Sh}_{p_1,p_2}^{p_1+p_2}}^{} \epsilon(\sigma) \, f_1(y_1,\dots,y_{p_1})f_2(z_1,\dots,z_{p_2}), \\
(f_1 \prec f_2)(x_1,\dots,x_{p_1+p_2}) &= \sum_{\sigma \in \mathrm{Sh}_{p_1,p_2}^{p_1}}^{} \epsilon(\sigma) \, f_1(y_1,\dots,y_{p_1})f_2(z_1,\dots,z_{p_2}),
\end{align*}
where 
$$
\left\{\begin{array}{l} y_k = x_{i_1} \rhd \dots \rhd x_{i_j} \rhd x_{\sigma(k)}, \\
z_k = x_{\sigma(p_1+k)},
\end{array}
\right.
$$
with $x_{i_1} < \dots < x_{i_j} < \sigma(k)$ and $i_l \in \{\sigma(p_1+1),\dots,\sigma(p_1+p_2)\}$,
is a differential graded dendriform algebra. In particular, $HR^{\bullet}(X,A)$ is provided with a graded dendriform algebra structure.
\end{theorem}
\begin{remark}
If the rack is pointed, then these formulas are well defined on the subcomplex $\{CR^n_p(X,A),d_R^n\}_{n \in \mathbb{N}}$. Notice that these formulas hold also in the case of a shelf (resp. pointed shelf).
\end{remark}
\section{A graded associative algebra morphism from $H^{\bullet}(G,A)$ to $HR^{\bullet}(\mathrm{Conj}(G),A)$}
Let $G$ be group and $A$ be an associative algebra over $\mathbb{Z}$ considered as a trivial $G$-module. The cochain complex $\{C^n(G,A),d_G^n,\cup\}$ calculating the group cohomology with coefficients in $A$ is provided with a differential graded associative algebra structure given by the cup product $\cup$. Moreover, considering the rack $\mathrm{Conj(G)}$ associated to $G$ (Example \ref{Example : groups and augmented racks}), we have shown (Theorem \ref{Theorem : CR(X,A) is a differential graded dendriform algebra}) that the cochain complex $\{CR^n(\mathrm{Conj(G)},A),d_R^n,\star\}$, with $\star = \, \succ + \prec$, is a differential graded associative algebra. In this section, we define a differential graded associative algebra morphism from $\{C^n(G,A),d_G^n,\cup\}$ to $\{CR^n(\mathrm{Conj(G)},A),d^n_R,\star\}$.
\paragraph{Group cohomology.} Let $G$ be a group and $A$ be a left $G$-module. The cochain complex $\{C^n(G,A),d_G^n\}_{n \in \mathbb{N}}$ is defined by
\begin{align*}\displaystyle
C^n(G,A) &:= \mathrm{Hom}_{\mathbf{Set}}(G^n,A),\\
d_G^{n+1} &:= \sum_{i=0}^{n+1} (-1)^i d_i^{n+1}f,
\end{align*}
where
$$
d_i^{n+1}f(x_1,\dots,x_{n+1}) := \left\{\begin{array}{ll}
x_1\cdot f(x_2,\dots,x_{n+1}) & \textrm{if } i = 0,\\
f(x_1,\dots,x_ix_{i+1},\dots,x_{n+1}) & \textrm{if } 1 \leq i \leq n,\\
f(x_1,\dots,x_{n}) & \textrm{if } i = n+1.
\end{array}
\right.
$$
The family of maps $\{d_i^n\}$ satisfied the \textit{simplicial identities}, that is, for all $0 \leq i < j \leq n+1$, we have the identities :
$$
d_i^{n+1} \circ d_{j-1}^n = d_j^{n+1} \circ d_{i}^n,
$$
so $d_G$ is a differential.
\par
The cohomology associated to this cochain complex is called the \textit{group cohomology of $G$ with coefficients in $A$} and denoted $H^{\bullet}(G,A)$.
\paragraph{Cup product on group cohomology.} Let $G$ be a group and $A$ be an associative algebra with product denoted by $\mu_A$. Let $p_1,p_2 \in \mathbb{N}$, let $\rho$ denote the set map from $G^{p_1+p_2}$ to $G^{p_1} \times G^{p_2}$ defined by
$$
\rho(x_1,\dots,x_{p_1+p_2}) := \big((x_1,\dots,x_{p_1}),(x_{p_1+1},\dots,x_{p_1+p_2})\big).
$$
A graded product $\cup$, called \textit{cup product}, is defined on $\{C^n(G,A),d_G^n\}_{n \in \mathbb{N}}$ by the formula
$$\
f_1 \cup f_2 := \mu_A \circ (f_1 \times f_2) \circ \rho. 
$$
\begin{theorem}
Let $G$ be a group and $A$ be an associative algebra (considered as a trivial $G$-module). The cochain complex $\{C^n(G,A),d_G^n\}_{n \in \mathbb{N}}$ provided with the cup product $\cup$ is a differential graded associative algebra structure. 
\end{theorem}
\paragraph{Equivalent definitions of the group cohomology and the cup product.} Let $G_{\mathbf{Cat}}$ denote the category canonically associated to a group $G$ (Example \ref{Example : category structure on square_n}). Let $\mathrm{B}(G)$ denote its \textit{presimplicial nerve of $G_{\mathbf{Cat}}$}, that is, the functor from $\Delta^{op}$ to $\mathbf{Set}$ defined by
\begin{align*}
\mathrm{B}(G)(\Delta_n) &:= \mathrm{Hom}_{\mathbf{Cat}}(\Delta_n,G_{\mathbf{Cat}}),\\
\mathrm{B}(G)(\partial_i^{n}) &:=  \mathrm{Hom}_{\mathbf{Cat}}(-,G_{\mathbf{Cat}})(\partial_i^n) = (\partial_i^{n})^*.
\end{align*}
\begin{proposition}
Let $G$ be a group. There is a bijection $\mathrm{B}(G)(\Delta_n) \stackrel{\lambda_n}{\simeq} G^n$, and under this bijection we have
$$
(\partial_i^{n})^*(x_1,\dots,x_n) = \left\{\begin{array}{ll}
(x_2,\dots,x_{n+1}) & \textrm{if } i = 0,\\
(x_1,\dots,x_ix_{i+1},\dots,x_{n+1}) & \textrm{if } 1 \leq i \leq n,\\
(x_1,\dots,x_{n}) & \textrm{if } i = n+1.
\end{array}
\right.
$$ 
\end{proposition}
Let $A$ be a left $G$-module. The composition $\mathrm{Hom}_{\mathbf{Set}}(-,A) \circ \mathrm{B}(G)$ defines a precosimplicial abelian group . Denote by $\{C^n_{\Delta}(G,A),d_{\Delta}^n\}_{n \in \mathbb{N}}$ the cochain complex associated to this precosimplicial abelian group. The cohomology associated to this cochain complex is denoted $H^{\bullet}_{\Delta}(G,A)$.
\begin{corollary}\label{Corollary : isomorphism between C(G,A) and CDelta(G,A)}
Let $G$ be a group and $A$ be an abelian group (considered as a trivial left $G$-module). There is an isomorphism of cochain complexes
$$
\{C^n(G,A),d_G^n\}_{n \in \mathbb{N}} \stackrel{\lambda^{\bullet}}{\simeq} \{C^n_{\Delta}(G,A),d_{\Delta}^n\}_{n \in \mathbb{N}}.
$$
\end{corollary}
Let $p_1,p_2 \in \mathbb{N}$. We define two functors $j_{p_1} : \Delta_{p_1} \to \Delta_{p_1+p_2}$ and $j_{p_2} : \Delta_{p_2} \to \Delta_{p_1+p_2}$ by $j_{p_1}(k) = k$ and $j_{p_2}(k) = p_1+k$. Let $\rho$ denote the set map from $\mathrm{B}(G)(\Delta_{p_1+p_2})$ to $\mathrm{B}(G)(\Delta_{p_1}) \times \mathrm{B}(G)(\Delta_{p_2})$ defined by
$$
\rho := (j_{p_1}^*,j_{p_2}^*).
$$
A graded product $\cup_{\Delta}$ on $\{C^n_{\Delta}(G,A),d_{\Delta}^n\}_{n \in \mathbb{N}}$ is defined by
$$
f_1 \cup_{\Delta} f_2 := \mu_A \circ (f_1 \times f_2) \circ \rho.
$$
\begin{proposition}
Let $G$ be a group and $A$ be an associative algebra (considered as a trivial $G$-module). The cochain complex $\{C^n_{\Delta}(G,A),d_{\Delta}^n\}_{n \in \mathbb{N}}$ provided with the product $\cup_{\Delta}$ is a differential graded associative algebra.
\end{proposition}
\begin{proposition}
Let $G$ be a group and $A$ be an associative algebra (considered as a trivial $G$-module). The cochain complex isomorphism between $\{C^n(G,A),d_G^n\}_{n \in \mathbb{N}}$ and $\{C^n_{\Delta}(G,A),d_{\Delta}^n\}_{n \in \mathbb{N}}$ given in Corollary \ref{Corollary : isomorphism between C(G,A) and CDelta(G,A)} is a differential graded associative algebra isomorphism.
\end{proposition}
\paragraph{Cubical cohomology of groups.} We have seen in Example \ref{Example : category structure on square_n} and \ref{Example : functors} that $\square_n$ has a category structure and that $\partial_{i,\epsilon}^n$ is a functor. Hence, there is a precubical set $N(\mathcal{C}) : \square^{op} \to \mathbf{Set}$ associated to any category $\mathcal{C}$, called \textit{the cubical nerve of $\mathcal{C}$}, and defined by
\begin{align*}
N(\mathcal{C})(\square_n) &:= \mathrm{Hom}_{\mathbf{Cat}}(\square_n,\mathcal{C}),\\
N(\mathcal{C})(\partial_{i,\epsilon}^n) &:= \mathrm{Hom}_{\mathbf{Cat}}(-,\mathcal{C})(\partial_{i,\epsilon}^n) = (\partial_{i,\epsilon}^n)^*.
\end{align*}
Let $G$ be a group and $A$ be a left $G$-module. The composition $\mathrm{Hom}_{\mathbf{Set}}(-,A) \circ N(G_{\mathbf{Cat}})$ defines a precocubical abelian group. Denote by $\{C_{\square}^n(G,A),d_{\square}^n\}_{n \in \mathbb{N}}$ the cochain complex associated to this precocubical abelian group. The cohomology associated to this cochain complex is called the \textit{cubical cohomology of $G$ with coefficient in $A$} and denoted $H^{\bullet}_{\square}(G,A)$.
\paragraph{The morphism $H^{\bullet}(G,A) \stackrel{[S^{\bullet}]}{\rightarrow} HR^{\bullet}(G,A)$.} Define a graded abelian group morphism $S^{\bullet}$ as the following composition :

\begin{align}\label{definition of S}
\xymatrix{
C^{\bullet}(G,A) \ar[d]_{\lambda^{\bullet}} \ar@{.>}[rrr]^{S^{\bullet}} & & & CR^{\bullet}(\mathrm{Conj}(G),A)
\\
C_{\Delta}^{\bullet}(G,A) \ar[r]^{\Sigma^{\bullet}} & C_{\square}^{\bullet}(G,A) \ar[r]^-{T^{\bullet}} & C^{\bullet}(N(\mathrm{Tr}\big(G_{\mathbf{Cat}})\big),A) \ar[r]^{I^{\bullet}} & C^{\bullet}(N(\mathrm{Conj}(G)),A) \ar[u]^{(\eta^{\bullet})^{-1}}
} 
\end{align}
where
\begin{itemize} 
\item $I^{n} = \big(\mathrm{Hom}_{\mathbf{Set}}(-,A) \circ \mathrm{Hom}_{\mathbf{Trunk}}(\square_n,-)\big)(\mathrm{inc})$ with $\mathrm{inc} \in \mathrm{Hom}_{\mathbf{Trunk}}\big(\mathrm{Conj(G)},\mathrm{Tr}(G_{\mathbf{Cat}})\big)$ is the inclusion of $\mathrm{Conj(G)}$ into $\mathrm{Tr}(G_{\mathbf{Cat}})$,
\item $T^{n} = \mathrm{Hom}_{\mathbf{Set}}(-,A)\big(\theta(\square_{n},G_{\mathbf{Cat}})\big)$ with $\theta(\square_n,G_{\mathbf{Cat}})$ the bijection induced by the adjunction $\mathrm{Cat} \stackrel{\theta}{\vdash} \mathrm{Tr}$ (cf. \eqref{adjunction Cat dashv Tr}) and the equality $\mathrm{Cat}(\square_n) = \square_n$ (cf. \eqref {Cat(square)=square}).
\item $\displaystyle \Sigma^{n} = \sum_{\sigma \in \mathbb{S}_n}^{} \epsilon(\sigma) \, s_{\sigma}$ with $s_{\sigma} = \big(\mathrm{Hom}_{\mathbf{Set}}(-,A) \circ \mathrm{Hom}_{\mathbf{Cat}}(-,G_{\mathbf{Cat}})\big)(\sigma)$ where $\sigma \in \mathrm{Hom}_{\mathbf{Cat}}(\Delta_n,\square_n)$ is the functor defined in Example \eqref{Example : functors}. 
\end{itemize}
\begin{theorem}\label{Theorem : construction of the morphism S}
Let $G$ be a group and $A$ be an abelian group (considered as a trivial $G$-module). The graded linear map $S^{\bullet}$ is a chain complex morphism.
\end{theorem}
\begin{proof}
The map $I^{\bullet}$ and $T^{\bullet}$ are induced by cubical set morphisms, hence they are chain complex morphisms. It remains to prove that $\Sigma^{\bullet}$ is a chain complex morphism.
\par
We want to prove that $d^{n+1}_{\square}\big(\Sigma^n(f)\big) = \Sigma^{n+1}\big(d_{\Delta}^{n+1}(f)\big)$. The left hand term of this equation is equal to the sum of two terms (A) and (B) with
\begin{align*} \displaystyle
\mathrm{(A)} &= \sum_{i = 1}^{n+1} \sum_{\sigma \in \mathbb{S}_n}^{} (-1)^{i+1} \epsilon(\sigma) \, f \circ (\partial_{i,1}^{n+1} \circ \sigma)^*,\\
\mathrm{(B)} &= \sum_{i = 1}^{n+1} \sum_{\sigma \in \mathbb{S}_n}^{} (-1)^i \epsilon(\sigma) \, f \circ (\partial_{i,0}^{n+1} \circ \sigma)^*.
\end{align*} 
The right hand term is equal to sum of three terms (a), (b) and (c) with
\begin{align*}\displaystyle
\mathrm{(a)} &= \sum_{\sigma \in \mathbb{S}_{n+1}}^{} \epsilon(\sigma) \, f \circ (\sigma \circ \partial_0^{n+1}),\\
\mathrm{(b)} &= \sum_{\sigma \in \mathbb{S}_{n+1}}^{} (-1)^{n+1} \epsilon(\sigma) \, f \circ (\sigma \circ \partial_{n+1}^{n+1}),\\
\mathrm{(c)} &= \sum_{\sigma \in \mathbb{S}_{n+1}}^{} \sum_{i = 1}^n (-1)^{i} \epsilon(\sigma) \, f \circ (\sigma \circ \partial_{i}^{n+1}).
\end{align*}
In order to prove the wanted equality we are going to show that
$$
\mathrm{(a)} = \mathrm{(A)}, \, \mathrm{(b)} = \mathrm{(B)} \, \text{ and } \, \mathrm{(c)} = 0.
$$
$\bullet \mathrm{(a)} = \mathrm{(A)}:$ This equality is proved using the change of variables
$$
\nu : \mathbb{S}_n \times \{1,\dots,n+1\} \longrightarrow \mathbb{S}_{n+1}
$$
defined by $\nu(\sigma,i) := (i \, i-1 \dots 1) \circ (1 \star \sigma)$. Using this change of variables, we have to prove for all $\sigma \in \mathbb{S}_n$ and $i \in \{1,\dots,n+1\}$ :
\begin{enumerate}
\item $(\partial_{i,1}^{n+1} \circ \sigma)^* = (\nu(\sigma,i) \circ \partial_{0}^{n+1})^*$,
\item $(-1)^{i+1} \epsilon(\sigma) = \epsilon\big(\nu(\sigma,i)\big)$.
\end{enumerate}
$\mathrm{1.}$ Let $\sigma \in \mathbb{S}_n$ and $i \in \{1,\dots,n+1\}$. We have $(\partial_{i,1}^{n+1} \circ \sigma)(k) = \sigma([k])_{< i} \amalg t_{+1}(\sigma([k])_{\geq i}) \amalg \{i\}$ and 
\begin{align*}
(\nu(\sigma,i) \circ \partial_0^{n+1})(k) &= \nu(\sigma,i)(k+1),\\
&= (i \, i-1 \dots 1)(\{1,1+\sigma(1),\dots,1+\sigma(k)\},\\
&= \{i\} \amalg t_{-1}\big((1+\sigma([k]))_{<i}\big) \amalg t_{+1}(\sigma([k]))_{\geq i}),\\
&= \sigma([k])_{< i} \amalg t_{+1}(\sigma([k])_{\geq i}) \amalg \{i\}.
\end{align*}
$\mathrm{2.}$ Clear.
\vskip 0.2 cm
\noindent
$\bullet \mathrm{(b)} = \mathrm{(B)}:$ This equality is proved using the change of variables
$$
\xi : \mathbb{S}_n \times \{1,\dots,n+1\} \longrightarrow \mathbb{S}_{n+1}
$$
defined by $\xi(\sigma,i) := (i \, i+1 \dots n+1) \circ (\sigma \star 1)$. Using this change of variables, we have to prove for all $\sigma \in \mathbb{S}_n$ and $i \in \{1,\dots,n+1\}$ :
\begin{enumerate}
\item $(\partial_{i,0}^{n+1} \circ \sigma)^* = (\xi(\sigma,i) \circ \partial_{n+1}^{n+1})^*$,
\item $(-1)^{i} \epsilon(\sigma) = (-1)^{n+1}\epsilon\big(\xi(\sigma,i)\big)$.
\end{enumerate}
$\mathrm{1.}$ Let $\sigma \in \mathbb{S}_n$ and $i \in \{1,\dots,n+1\}$. We have $(\partial_{i,0}^{n+1} \circ \sigma)(k) = \sigma([k])_{< i} \amalg t_{+1}(\sigma([k])_{\geq i})$ and 
\begin{align*}
(\xi(\sigma,i) \circ \partial_{n+1}^{n+1})(k) &= \xi(\sigma,i)(k),\\
&= (i \, i+1 \dots n+1)(\{\sigma(1),\dots,\sigma(k)\},\\
&= \sigma([k])_{<i} \amalg t_{+1}(\sigma([k]))_{\geq i}).
\end{align*}
$\mathrm{2.}$ Clear.
\vskip 0.2 cm
\noindent
$\bullet \mathrm{(c)} = 0:$ This equality is proved using the change of variables
$$
\kappa : \{\sigma \in \mathbb{S}_{n+1} \, | \, \epsilon(\sigma)=1\} \times \{1,\dots,n\} \longrightarrow \{\sigma \in \mathbb{S}_{n+1} \, | \, \epsilon(\sigma)=-1\} \times \{1 \dots,n\}
$$
defined by $\kappa(\sigma,i) := \big((\sigma(i) \, \sigma(i+1)\big) \circ \sigma,i)$. Using this change of variables, we have to prove for all $\sigma \in \mathbb{S}_{n+1}, \, \epsilon(\sigma)=1$ and $i \in \{1,\dots,n\}$:
\begin{enumerate}
\item $(\sigma \circ \partial_{i}^{n+1})^* = (\kappa(\sigma,i) \circ \partial_{i}^{n+1})^*$,
\item $(-1)^i \epsilon(\sigma) + (-1)^{i}\epsilon\big(\kappa(\sigma,i)\big)$.
\end{enumerate}
$\mathrm{1.}$ Let $\sigma \in \mathbb{S}_{n+1}$ and $i \in \{1,\dots,n\}$. We have 
\begin{align*}
(\kappa(\sigma,i) \circ \partial_{i}^{n+1})(k) &= \left\{ \begin{array}{ll}
									\kappa(\sigma,i)(k) & \textrm{if } k < i,\\
									\kappa(\sigma,i)(k+1) & \textrm{if } k \geq i.
									\end{array}
									\right.
								&= \left\{ \begin{array}{ll}
									\sigma([k]) & \textrm{if } k < i,\\
									\sigma([k+1]) & \textrm{if } k \geq i.
									\end{array}
									\right.
								&= (\sigma \circ \partial_{i}^{n+1})(k). 
\end{align*}
$\mathrm{2.}$ Clear.
\end{proof}
\paragraph{Explicit formula for $S^{\bullet}$.} By definition, for all $f \in C^n(G,A)$
$$\displaystyle
S^n(f) = \sum_{\sigma \in \mathbb{S}_n}^{} \epsilon(\sigma) \, f \circ  \lambda  \circ \sigma^* \circ \theta \circ \mathrm{inc}_* \circ \eta^{-1}.
$$
Then to have an explict formula for $S^{\bullet}$, we have to compute
$$
\lambda  \circ \sigma^* \circ \theta \circ \mathrm{inc}_* \circ \eta^{-1} : \mathrm{Conj}(G)^n \to G^n
$$
Let $(x_1,\dots,x_n) \in \mathrm{Conj}(G)^n$,
$$
(\lambda  \circ \sigma^* \circ \theta \circ \mathrm{inc}_* \circ \eta^{-1})(x_1,\dots,x_n) = (y_1,\dots,y_n),
$$
where $y_k = (F \circ \sigma)(k-1 \to k)$ for all $1 \leq k \leq n$ with $F = (\eta^{-1})(x_1,\dots,x_n)$. Thus using Lemma \ref{Lemma : nerve of a rack} we find
\begin{align*}
y_k &= (F \circ \sigma)(k-1 \to k),\\
&= F\big(\sigma([k-1]) \to \sigma([k])\big), \\ 
&= \displaystyle \prod_{\stackrel{1 \leq x \leq \sigma(k)}{x \notin \sigma([k-1])}} F([x-1] \to [x]), \\
&= x_{i_1} \rhd \dots \rhd x_{i_j} \rhd x_{\sigma(k)},
\end{align*}
with $x_{i_1} < \dots < x_{i_j} < \sigma(k)$ and $i_l \in \{\sigma(k+1),\dots,\sigma(n)\}$.
\begin{theorem}
Let $G$ be a group and $A$ be an abelian group (considered as a trivial $G$-module). The graded abelian group morphism $S^{\bullet}$ from $C^{\bullet}(G,A)$ to $CR^{\bullet}(\mathrm{Conj}(G),A)$ defined by
$$\displaystyle
S^n(f)(x_1,\dots,x_n) := \sum_{\sigma \in \mathbb{S}_n}^{} \epsilon(\sigma) \, f(y_1,\dots,y_n)
$$ 
with $y_k = x_{i_1} \rhd \dots \rhd x_{i_j} \rhd x_{\sigma(k)}$ where $x_{i_1} < \dots < x_{i_j} < \sigma(k)$ and $i_l \in \{\sigma(k+1),\dots,\sigma(n)\}$, is a cochain complex morphism.
\end{theorem}
\paragraph{$S^{\bullet}$ is a differential graded associative algebra morphism.} Let $G$ be a group and $A$ be an associative algebra (considered as a trivial $G$-module). There are two differential graded associative algebras :
\begin{itemize}
\item $\{C^n(G,A),d_G^n,\cup\}_{n \in \mathbb{N}}$,
\item $\{CR^n(\mathrm{Conj}(G),A),d_R^n,\star = \, \succ + \prec\}_{n \in \mathbb{N}}$.
\end{itemize}
and a cochain complex morphism $\{S^n\}_{n \in \mathbb{N}}$ between these two cochain complexes is defined in \eqref{definition of S}. The following theorem states that $\{S^n\}_{n \in \mathbb{N}}$ respects their associative structures. 
\begin{theorem}\label{Theorem : a graded associative algebra morphism from group cohomology to rack cohomology}
Let $G$ be a group and $A$ be an associative algebra (considered as a trivial $G$-module). The cochain complex morphism 
$$
\{S^{n}\}_{n \in \mathbb{N}} : \{C^{n}(G,A),d^n_G\}_{n \in \mathbb{N}} \longrightarrow \{CR^{n}(\mathrm{Conj}(G),A),d^n_R\}_{n \in \mathbb{N}}
$$ 
is a differential graded associative algebra morphism. Hence it induces a graded associative algebra morphism $\{[S^{n}]\}_{n \in \mathbb{N}}$ from $\{H^{n}(G,A),\cup\}_{n \in \mathbb{N}}$ to $\{HR^n(\mathrm{Conj}(G),A),\star\}_{n \in \mathbb{N}}$.
\end{theorem}
\begin{proof}
By construction $\lambda^{\bullet}$ and $\eta^{\bullet}$ are graded associative algebra morphisms, hence we have to prove that $I^{\bullet} \circ T^{\bullet} \circ \Sigma^{\bullet}$ is a graded associative algebra morphism.
\par
Let $f_1 \in C_{\Delta}^{p_1}(G,A)$ and $f_{2} \in C_{\Delta}^{p_2}(G,A)$, on one hand 
$$\displaystyle
(I^{p_1+p_2} \circ T^{p_1+p_2} \circ \Sigma^{p_1+p_2})(f_1 \cup_{\Delta} f_2) := \sum_{\sigma \in \mathbb{S}_{p_1+p_2}}^{} \epsilon(\sigma) \, \mu_A \circ (f_1 \times f_2) \circ \rho \circ \sigma^* \circ \theta \circ \mathrm{inc}_*, 
$$
and on the other hand
\begin{align*}\displaystyle
(I^{p_1} \circ T^{p_1} \circ \Sigma^{p_1})&(f_1) \star (I^{p_2} \circ T^{p_2} \circ \Sigma^{p_2})(f_2) := \\  
& \sum_{\gamma \in \mathrm{Sh}_{p_1,p_2}}^{} \sum_{\alpha \in \mathbb{S}_{p_1}}^{} \sum_{\beta \in \mathbb{S}_{p_2}} \epsilon(\gamma) \epsilon(\alpha) \epsilon(\beta) \, \mu_A \circ (f_1 \times f_2) \circ \big((\alpha^* \circ \theta \circ \mathrm{inc}_*) \times (\beta^* \circ \theta \circ \mathrm{inc}_*)\big) \circ \rho_{\sigma}.
\end{align*}
To show the equality of these two sums we are going to use the change of variables
$$
\varphi : \mathrm{Sh}_{p_1,p_2} \times \mathbb{S}_{p_1} \times \mathbb{S}_{p_2} \longrightarrow \mathbb{S}_{p_1+p_2}
$$
defined by $\varphi(\gamma,\alpha,\beta) = \gamma \circ (\alpha \star \beta)$. Using this change of variables we have to prove for all $\gamma \in \mathrm{Sh}_{p_1,p_2}$, $\alpha \in \mathbb{S}_{p_1}$ and $\beta \in \mathbb{S}_{p_2}$ :
\begin{enumerate}
\item $\rho \circ \sigma^* \circ \theta \circ \mathrm{inc}_* = \big((\alpha^* \circ \theta \circ \mathrm{inc}_*) \times (\beta^* \circ \theta \circ \mathrm{inc}_*)\big) \circ \rho_{\sigma}$,
\item $\epsilon(\gamma)\epsilon(\alpha)\epsilon(\beta) = \epsilon\big(\gamma \circ (\alpha \star \beta)\big)$.
\end{enumerate}
$\mathrm{1.}$ Let $\gamma \in \mathrm{Sh}_{p_1,p_2}$, $\alpha \in \mathbb{S}_{p_1}$ and $\beta \in \mathbb{S}_{p_2}$. We have to prove for all $F \in N\big(\mathrm{Conj}(G)\big)(\square_{p_1+p_2})$ the following equalities :
$$
\left\{\begin{array}{ll}
\big(F \circ \sigma \circ (\alpha \star \beta) \circ j_{p_1} \big)(k-1 \to k) = \big(F \circ \sigma \circ i_{p_1} \circ \alpha \big)(k-1 \to k)  & \textrm{if } 1 \leq k \leq p_1,\\
\big(F \circ \sigma \circ (\alpha \star \beta) \circ j_{p_2} \big)(k-1 \to k) = \big(F \circ \sigma \circ i_{p_2} \circ \beta \big)(k-1 \to k)  & \textrm{if } 1 \leq k \leq p_2.\\
\end{array}
\right.
$$
Let $1 \leq k \leq p_1$,
$$
\big(\sigma \circ (\alpha \star \beta) \circ j_{p_1} \big)(k) = (\sigma \circ (\alpha \star \beta)\big)(k) = (\sigma \circ (\alpha \star \beta)\big)([k]) =  \sigma\big(\alpha([k])\big)
$$ 
and 
$$
(\sigma \circ i_{p_1} \circ \alpha)(k) = (\sigma \circ i_{p_1})\big(\alpha([k])\big) = \sigma\big(\alpha([k])\big).
$$
Let $1 \leq k \leq p_2$, 
$$
\big(\sigma \circ (\alpha \star \beta) \circ j_{p_2} \big)(k) = (\sigma \circ (\alpha \star \beta)\big)(p_1+k) = (\sigma \circ (\alpha \star \beta)\big)([p_1+k]) =  \sigma\big(t_{p_1}\beta([k])\big)
$$ 
and 
$$
(\sigma \circ i_{p_2} \circ \beta)(k) = (\sigma \circ i_{p_2})\big(\beta([k])\big) = \sigma\big(t_{p_1}\beta([k])\big).
$$
$\mathrm{2.}$ Clear.
\end{proof}
\bibliographystyle{alpha}
\bibliography{bibliographie}
\end{document}